\documentclass[12pt]{amsart}

\usepackage{amsmath,amsthm,verbatim,amssymb,amsfonts,amscd, graphicx,cite,bbm,url,enumerate,hyperref,comment,xspace,mathtools}
\usepackage{graphics}
\hypersetup{
	colorlinks=true,
	linkcolor=blue,
	citecolor=black
}
\theoremstyle{plain}
\newtheorem{theorem}{Theorem}[section]

\newtheorem{theoremnb}{Theorem}
\newtheorem{lemma}[theorem]{Lemma}
\newtheorem{proposition}[theorem]{Proposition}

\theoremstyle{definition}

\theoremstyle{definition}

\allowdisplaybreaks[2]
\numberwithin{equation}{section}

\newcommand{\pt}{\, \dashv \,}

\begin{document}

\begin{abstract}
    The goal of this work is to prove an analogue of a recent result of Harper on almost sure lower bounds of random multiplicative functions, in a setting that can be thought of as a simplified function field analogue. It answers a question raised in work of Soundararajan and Zaman, who proved moment bounds for the same quantity in analogy to those of Harper in the random multiplicative setting. Having a simpler quantity allows us to make the proof close to self-contained, and perhaps somewhat more accessible.
\end{abstract}

\title[Almost sure lower bounds for a model problem for multiplicative chaos]{Almost sure lower bounds for a model problem for multiplicative chaos in number theory}

\author{Maxim Gerspach}
\address{Department of Mathematics, KTH Royal Institute of Technology, Lindstedtsvägen 25, 114 28 Stockholm, Sweden}
\email{gerspach@kth.se}

\thanks{This work received funding from the Göran Gustafsson Foundation and from the Swedish Research Council Grant No. 2016-05198.}

\date{\today}

\maketitle

\section{Introduction}

Let $(X(k))_{k \ge 1}$ be a sequence of independent standard complex Gaussian random variables, i.e. so that their real and imaginary parts are independent, centered and have variance $\frac{1}{2}$. Consider the sequence of random variables $(A(n))_{n \ge 0}$ defined through the following formal identity of power series
\[ \exp \Big( \sum_{k \ge 1} \frac{X(k)}{\sqrt{k}} z^k \Big) = \sum_{n \ge 0} A(n) z^n. \]
These random variables have recently been considered as an analogue of a (Steinhaus) random multiplicative function that seems to be rich enough to share most of their properties, while often allowing for simpler proofs thereof. Recall that a Steinhaus random multiplicative function is a (random) completely multiplicative function $f : \mathbb{N} \to \mathbb{C}$ such that its values at the primes are independent and uniformly distributed on the complex unit circle (Steinhaus distributed).

In this setting, there has been a lot of attention in recent years regarding the moments and almost sure bounds for mean values of random multiplicative functions. In particular, Harper proved in \cite{Harper8} that for any function $V(x)$ tending to infinity with $x$ there are, almost surely, arbitrarily large values of $x$ such that
\begin{equation}\label{HarperASLB} \Big| \frac{1}{\sqrt{x}} \sum_{n \le x} f(n) \Big| \ge \frac{(\log \log x)^{1/4}}{V(x)}. \end{equation}
Regarding upper bounds, the best known result to this date is due to Lau, Tenenbaum and Wu \cite{LTW1} and states that, almost surely, we have
\[ \frac{1}{\sqrt{x}} \sum_{n \le x} f(n) = O_\varepsilon((\log \log x)^{2 + \varepsilon}). \]
Moreover, Mastrostefano \cite{Mastrostefano1} recently proved an upper bound essentially matching the one of Harper for the same sum when restricted to integers that possess a large prime factor. More precisely, denoting by $P(n)$ the largest prime factor of an integer $n$, he showed that we almost surely have
\[ \frac{1}{\sqrt{x}} \sum_{\substack{ n \le x \\ P(n) > \sqrt{x} }} f(n) = O_\varepsilon((\log \log x)^{1/4 + \varepsilon}). \]

These almost sure bounds may be compared to the first moment: It is known, as a corollary to \cite[Theorem 1]{Harper1}, that we have
\[ \mathbb{E} \bigg[ \, \bigg| \frac{1}{\sqrt{x}} \sum_{n \le x} f(n) \bigg | \bigg] \asymp \frac{1}{(\log \log x)^{1/4}}. \]
This (suspected) discrepancy of a factor $\sqrt{\log \log x}$ between the first moment and almost sure behaviour is akin to the law of the iterated logarithm for independent random variables. We refer the reader to the discussion in \cite{Harper8} for more details on this comparison.

Inspired by this result, Soundararajan and Zaman consider in \cite{SoundZaman1} the moments of $A(n)$, and showed that they behave rather similar to the random multiplicative setting. They proved in particular that we have
\[ \mathbb{E}[|A(n)|] \asymp \frac{1}{(\log n)^{1/4}}. \]
They further raise the question whether one can obtain analogues of the almost sure results on random multiplicative functions for $A(n)$.

The main goal of this work is to prove that this is indeed the case, i.e. to prove the following analogue of \eqref{HarperASLB}.

\begin{theoremnb}\label{T1}
    For any function $V(n)$ tending to infinity with $n$, there almost surely exist arbitrarily large values of $n$ for which
    \[ |A(n)| \ge \frac{(\log n)^{1/4}}{V(n)}. \]
\end{theoremnb}

The link between mean values of random multiplicative functions and the random variables $A(n)$ is perhaps not obvious at first sight. As explained by Soundararajan and Zaman in \cite{SoundZaman1}, it is more apparent in a function field setting of $\mathbb{F}_q[t]$, especially in the limit as $q \to \infty$. This analogy roughly corresponds to imposing that all primes $e^k < p \le e^{k+1}$ shall have ``size'' $e^k$, which is more sensible in a function field setting, where one tends to have many irreducible polynomials of the same degree. We refer the reader to their introduction for the details of this analogy. In addition, we will give a more detailed heuristic in the next section after introducing further definitions. 

As in the work of Harper \cite{Harper8}, we will prove this by showing a lower bound for the maximum of $|A(n)|$ on adequate intervals of $n$, to hold with high probability, and then deduce the Theorem from an invocation of the first Borel-Cantelli Lemma.

\begin{theoremnb}\label{T2}
    Uniformly for all large $N$ and all $1 \le W \le \frac{1}{10} \log \log N$, we have
    \[ \max_{\frac{8 N}{7} \le n \le \frac{4N}{3}} |A(n)| \ge \frac{(\log N)^{1/4}}{e^{6 W/5}} \]
    with probability $ \ge 1 - O(e^{-W/10})$.
\end{theoremnb}

In fact, the deduction of the main result from this Theorem follows \cite{Harper8} to the letter aside from some choices of constants, but will be included here for completeness.

\begin{proof}[Proof of Theorem \ref{T1} assuming Theorem \ref{T2}]
    Assume without loss of generality that $V(n) \le (\log n)^{1/10}$ for sufficiently large $n$. Set
    \[ W(N) := \min_{n \in [\frac{8N}{7}, \frac{4N}{3}]} \frac{5 \log V(n)}{6} - 1, \]
    so that $W(N)$ goes to infinity with $N$, but also satisfies $W(N) \le \frac{1}{10} \log \log N$ when $N$ is sufficiently large. Hence, by Theorem \ref{T2}, the probability that for a given sufficiently large $N$ there is no $n \in [\frac{8N}{7}, \frac{4N}{3}]$ such that \[|A(n)| \ge \frac{(\log N)^{1/4}}{e^{6 W(N)/5}} \ge \frac{(\log n)^{1/4}}{V(n)} \]
    is $O(e^{-W(N)/10})$ (since $n$ can also be assumed sufficiently large). Hence, by choosing an adequate subsequence of values of $N$ we can ensure these probabilities to be summable. Thus, the first Borel-Cantelli Lemma tells us that, almost sure, only finitely many of these complementary events occur on our subsequence. Excluding these finitely many events, we have (almost surely) found our sequence of $n$ going to infinity that satisfies the claimed inequality.
\end{proof}

Overall, our proof will closely follow the works of Harper \cite{Harper8, Harper1}, but there will be some simplifications arising from the fact that we are often dealing with Gaussian random variables in place of random variables that are approximately Gaussian. This will allow the argument to be essentially self-contained with somewhat less effort than in these works. We will keep many of the constants analogous to the ones that arise there, mostly for easier comparison, since we do have even more flexibility in their choices. We will moreover employ some of the ideas from the work of Soundararajan and Zaman \cite{SoundZaman1}.

One of the differences to the random multiplicative function setting is that there one needs to restrict to a (geometric) subsequence of values of $x$ in $[X^{8/7},X^{4/3}]$ in order to have any hope that $\frac{1}{\sqrt{x}}\sum_{n \le x} f(n)$ at different $x$ have small correlations. In our setting, this would correspond to taking a linear (arithmetic) subsequence, and it turns out that one can just as well use all values of $A(n)$ at this step.

It seems very much a possibility that one can also transfer the aforementioned almost sure upper bounds \cite{LTW1,Mastrostefano1} to the random variables $A(N)$ by going through the respective arguments and making similar adjustments as in this work or as in \cite{SoundZaman1}.

Finally, we also want to point the reader towards the work \cite{ASVZZ}, where the authors numerically study the asymptotics of the first moment of $A(n)$ (which are not known to exist) by developing an algorithm that allows a more effective computation.

\section{Preliminaries and basic estimates}

Let $t \in \mathbb{R}$ and $m \in \mathbb{N}$. Let $\sigma \ge 0 $, which the reader should think of as small (certainly $o(1)$). In the following, we will frequently be working with the random variables
\[ Z_t(m) = Z_{t, \sigma}(m) := \sum_{e^{m-1} < k \le e^m} \frac{\Re X(k) e^{i k t}}{\sqrt{k} e^{k \sigma}}, \]
where (somewhat arbitrarily) we will throughout include the term $k = 1$ in the sum when $m = 1$, and in corresponding sums. We also set $Z(m) := Z_0(m)$.
We note that these random variables have mean $0$ and variance
\[ \sigma_m^2 = \sigma_{m, \sigma}^2 := \mathbb{E}[Z_t(m)^2] = \sum_{e^{m-1} < k \le e^m} \frac{1}{2 k e^{2 k \sigma}}, \]
independent of $t$ (we hope that little confusion can arise from these two uses of the letter $\sigma$). Moreover, we define
\[ \rho_{m, t} \sigma_m^2 = \rho_{m, t, \sigma} \sigma_m^2 := \mathbb{E}[Z(m) Z_t(m)] = \sum_{e^{m - 1} < k \le e^m} \frac{\cos (k t)}{2 k e^{2 k \sigma}} \]
and remark that $\mathbb{E}[Z_t(m) Z_u(m)] = \rho_{m, u-t} \sigma_m^2$.

It will be useful to derive an explicit expression for $A(n)$ (which also appears in \cite{SoundZaman1}). Note that in the implicit definition, we can expand the exponential and see that
\begin{align*}
    \sum_{n \ge 0} A(n) z^n &= \sum_{m \ge 0} \frac{1}{m!} \Big( \sum_{k \ge 1} \frac{X(k)}{\sqrt{k}} z^k \Big)^m \\
    &= \sum_{m \ge 0} \frac{1}{m!} \sum_{n \ge 0} z^n \sum_{k_1 + \dots + k_m = n} \frac{X(k_1) \cdots X(k_m)}{\sqrt{k_1 \cdots k_m}} 
\end{align*}
Next, we can order the $k_i$ decreasingly and note that the number of ways of rearranging them is
\[ \frac{m!}{\prod_{j \ge 1} (\#\{k_i = j\})!  }. \]
In this way we have linked the previous expression to partitions of $n$, and simply comparing the coefficients in the respective power series we see that
\[ A(n) = \sum_{\lambda \pt n} \prod_{j \ge 1} \left( \frac{X(j)}{\sqrt{j}} \right)^{m_j} \frac{1}{m_j!} =: \sum_{\lambda \pt n} a(\lambda), \]
where we denote by $m_j = m_j(\lambda)$ the number of parts of $\lambda$ of size $j$.

Further, for $K > 1$ real and $z \in \mathbb{C}$, we set
\[ F_K(z) := \exp \bigg( \sum_{k \le K} \frac{X(k)}{\sqrt{k}} z^k \bigg).  \]
The reader should think of $z$ as being on or close to the complex unit circle. Most of the time we will be able to assume that $\log K$ is an integer, which implies in particular that
\[ |F_K(e^{-\sigma + i t})| = \exp \bigg( \sum_{m = 1}^{\log K} Z_t(m) \bigg). \]
Moreover, one verifies by the same computation as for $A(n)$ that we have
\[ F_K(z) = \sum_{n \ge 0} \Big( \sum_{\substack{\sigma \pt n \\ \sigma_1 \le K } } a(\sigma) \Big) z^n. \]
In particular, these coefficients are simply $A(n)$ whenever $n \le K$.

The way the analogy between the random variables $A(n)$ and the random multiplicative function setting manifests itself in the course of proof is through application of Cauchy's Theorem and Perron's formula, respectively. To explain this in more detail, note that for $N \le K$, Cauchy's formula gives
\begin{equation}\label{Cauchy} A(N) = \frac{1}{2 \pi i} \int_{|z| = 1} F_K(z) \frac{dz}{z^{N+1}}. \end{equation}

Let us compare this to
\[ \frac{1}{\sqrt{e^N}} \sum_{n \le e^N} f(n). \]
Define
\[ \tilde{F}_K(s) := \prod_{p \le e^K} \left( 1- \frac{f(p)}{p^s} \right)^{-1} = \sum_{P(n) \le e^K} \frac{f(n)}{n^s} \]
to be the (random and partial) Euler product associated to $f$.
Perron's formula tells us that
\[ \frac{1}{\sqrt{e^N}} \sum_{n \le e^N} f(n) \approx \frac{1}{2 \pi i e^{N/2}} \int_{(1)} \tilde{F}_K(s) e^{s N} \frac{ds}{s}, \]
where $(c)$ denotes the line with real part equal to $c$ (and in fact we have equality, but that is not important for this informal discussion). One then verifies (as we will later on) that this integral can be cut off at (say) $|\Im s| = e^{3 N/4}$, and that one can moreover shift the resulting integral to the line $\Re s = \frac{1}{2}$, with negligible error in both cases. Thus,
\[\frac{1}{\sqrt{e^N}} \sum_{n \le e^N} f(n) \approx \int_{|t|< e^{3N/4}} \tilde{F}_K \left( \frac{1}{2} + i t \right) e^{i N t} \frac{d t}{\frac{1}{2}+ i t}. \]
Ideally one would like to cut off this integral at a constant (such as $\pi$), and one can show that with high probability one can at least restrict to $|t| < N^2$ (and even $(\log N)^2$ in an appropriate sense). Comparing this to \ref{Cauchy} after setting $z = e^{-i t}$, these expressions seem to look rather similar if one chooses to believe that $F_K(e^{-it})$ and $\tilde{F}_K(1/2+it)$ behave in a similar fashion. And indeed, note that
\[ \log \tilde{F}_K(1/2 + i t) =  - \sum_{ p \le e^K} \log \left( 1 - \frac{f(p)}{p^{1/2+i t}} \right) \approx \sum_{p \le e^K} \frac{f(p)}{p^{1/2+i t}} = \sum_{k \le K} \sum_{e^{k-1} < p \le e^k} \frac{f(p)}{p^{1/2+i t}}. \]
For any $t$ and sufficiently large $k$, these sums over $p$ are a sum of many random variables with mean $0$ and variance $\frac{1}{p}$ (whose real and imaginary parts are uncorrelated and have variance $\frac{1}{2p}$), thus the central limit Theorem implies that for each $k$,
\[ \sum_{e^{k-1} < p \le e^k} \frac{f(p)}{p^{1/2+i t}} \]
is approximately a (complex) Gaussian with mean $0$ and variance
\[ \sum_{e^{k-1} < p \le e^k} \frac{1}{p} \approx \log k - \log(k - 1) \approx \frac{1}{k}  \]
(whose real and imaginary part are independent real Gaussians with mean $0$ and variance $\approx \frac{1}{2k})$.
The same holds true for $\frac{X(k)}{\sqrt{k}} e^{i t k}$ for any $t$. In fact one can also identify the correlations of the respective processes for different values of $t$.

We continue by recording the following elementary estimates.



\begin{proposition}\label{NTR2}
        For $\frac{1}{2} \le x < y$ (say), $\sigma \ge 0$ and $t \in (0,2\pi)$, we have
        \[ \bigg| \sum_{x < k \le y} \frac{e^{i t k}}{k e^{2 k \sigma}} \bigg| \le \frac{3 \pi}{x \Vert t \Vert}. \]
\end{proposition}

\begin{proof}
    Note that
    \[ \bigg| \sum_{k \le u} e^{i t k} \bigg| \le \frac{2}{|1-e^{it}|} \le \frac{\pi}{\Vert t \Vert} \]
    for all (say) $u \ge \frac{1}{2}$ and $t \in (0,2\pi)$. The claim then follows from partial summation.
\end{proof}

\begin{proposition}\label{RhoEstimate}
    Let $t \in \mathbb{R}$, $m \in \mathbb{N}$ and $\sigma \ge 0$. Then we have
    \[ e^{-2 \sigma e^m} \left( \frac{1}{2} + O(e^{-m}) \right) \le  \sigma_m^2 \le e^{-2 \sigma e^{m-1}} \left( \frac{1}{2} + O(e^{-m}) \right). \]
    If moreover $0 \le \sigma \le 10e^{-m}$ (say), then we have $\sigma_m^2 \asymp 1$ and
    \[ \rho_{m,t} \ll \frac{1}{\Vert t \Vert e^m}. \]
\end{proposition}

\begin{proof}
    Recall that we are still following the convention that the term corresponding to $k = 1$ is included in the sum when $m = 1$.
    
    To deduce the first part, one simply bounds $e^{-2 k \sigma}$ trivially and uses that
    \[ \sum_{e^{m-1} < k \le e^m} \frac{1}{2k} = \frac{1}{2} + O(e^{-m}). \]
    We note that this sum is in fact always between $\frac{1}{4}$ and $1$, say, as follows e.g. from \cite[Theorem]{RYoung} and a short computation. This certainly gives the next claim after bounding $e^{-2 k \sigma}$ trivially again.
    
    As for the estimate on $\rho_{m,t}$, note first that
    \[ \rho_{m,t} \asymp \rho_{m,t} \sigma_m^2 = \Re \sum_{e^{m-1} \le k < e^m} \frac{e^{i k t}}{2 k e^{2 k \sigma}}. \]
    One then applies partial summation, using again that
    \[ \sum_{k \le u} e^{i t k} \ll \frac{1}{\Vert t \Vert} \]
    and the claim follows.
\end{proof}

\begin{proposition}\label{MomentsZ}
    For any $m \in \mathbb{N}$, real $\alpha_1, \dots, \alpha_m$ and $t_1, \dots, t_m$, for $\frac{1}{2} \le x \le y$ (say) and $\sigma  \ge 0$, we have
    \begin{align*}
        \mathbb{E} \bigg[ \exp \bigg( \sum_{j = 1}^m 2 \alpha_j \sum_{x < k \le y} \frac{\Re X(k) e^{i k t_j}}{\sqrt{k} e^{i k \sigma}} \bigg) \bigg] = \exp \bigg( \sum_{j_1, j_2 = 1}^m \sum_{x < k \le y} \alpha_{j_1} \alpha_{j_2} \frac{\cos \big( k (t_{j_2}-t_{j_1}) \big)}{k e^{2 i k \sigma}} \bigg).
    \end{align*}
    In particular, uniformly for $\frac{1}{2} \le x \le y$, $\sigma \ge 0$ and $t_1, t_2 \in \mathbb{R}$, we have
    \begin{align*}
        \mathbb{E} \bigg[ \exp \bigg( \sum_{x < k \le y} \frac{\Re X(k) e^{i k t_1}}{\sqrt{k} e^{i k \sigma}} + \sum_{x < k \le y} \frac{\Re X(k) e^{i k t_2}}{\sqrt{k} e^{i k \sigma}} \bigg) \bigg] \ll \sqrt{\frac{y}{x}} \min \bigg \{ \sqrt{\frac{y}{x}}, \sqrt{\frac{1}{\Vert t_2 - t_1 \Vert x}} \bigg \}.
    \end{align*}
\end{proposition}

\begin{proof}
    The first part of this statement is completely elementary since we are dealing with true Gaussian random variables, and only included for easier comparison to the work of Harper \cite{Harper8}. Namely, one notes that the sum
    \[ \sum_{j = 1}^m 2 \alpha_j \sum_{x < k \le y} \frac{\Re X(k) e^{i k t_j}}{\sqrt{k} e^{i k \sigma}} \]
    inside the exponential is still Gaussian, with mean $0$ and variance
    \[ \sum_{x < k \le y} \mathbb{E} \bigg[ \bigg( 2 \sum_{j = 1}^m \alpha_j \frac{\Re X(k) e^{i k t_j}}{\sqrt{k} e^{i k \sigma}} \bigg)^2 \bigg] = 4 \sum_{j_1, j_2 = 1}^m \alpha_{j_1} \alpha_{j_2} \frac{\mathbb{E}\big[ (\Re X(k) e^{i k t_{j_1}} ) ( \Re X(k) e^{i k t_{j_2}} ) \big]}{k e^{2 i k \sigma}}. \]
    One then verifies that
    \[ \mathbb{E}\big[ (\Re X(k) e^{i k t_{j_1}} ) ( \Re X(k) e^{i k t_{j_2}} ) \big] = \frac{\cos \big( k (t_{j_2} - t_{j_1}) \big)}{2}, \]
    e.g. by using that $\Re X(k) e^{i k t} = (\Re X(k)) \cos(k t) - (\Im X(k)) \sin(k t)$ as well as $\mathbb{E}[\Re X(k) \Im X(k)] = 0$ and $\cos(x - y) = \cos x \cos y + \sin x \sin y$. Noting that if $G$ is a Gaussian with mean $0$ and variance $\sigma^2$ then $\mathbb{E}[\exp(G)] = \exp(\sigma^2/2)$ gives the first claim.
    
    As for the second part, note that the first part, with $\alpha_1 = \alpha_2 = \frac{1}{2}$, implies
    \[ \mathbb{E} \bigg[ \exp \bigg( \sum_{x < k \le y} \frac{\Re X(k) e^{i k t_1}}{\sqrt{k} e^{i k \sigma}} + \sum_{x < k \le y} \frac{\Re X(k) e^{i k t_2}}{\sqrt{k} e^{i k \sigma}} \bigg) \bigg] = \exp \bigg( \sum_{x < k \le y} \frac{1 + \cos \big( k (t_2 - t_1) \big)}{2 k e^{2 k \sigma}} \bigg). \] 
    Moreover, we certainly have
    \[ \sum_{ x < k \le y } \frac{1}{2 k e^{2 k \sigma}} \le \frac{1}{2}(\log y - \log x) + O(1), \]
    which gives rise to the factor $\sqrt{\frac{y}{x}}$ in our estimate and also the factor $\sqrt{\frac{y}{x}}$ inside the minimum, simply by bounding the cosine by $1$. We may thus assume from now on that $\Vert t_2 - t_1 \Vert \ge \frac{1}{y}$, and note that by Proposition \ref{NTR2} we have
    \[ \sum_{\frac{1}{\Vert t_2 - t_1 \Vert} < k \le y} \frac{\cos \big( k (t_2 - t_1) \big)}{2 k e^{2 k \sigma}} = \Re \sum_{\frac{1}{\Vert t_2 - t_1 \Vert} < k \le y} \frac{e^{i k (t_2 - t_1)}}{2 k e^{2 k \sigma}} \ll 1. \]
    If $\Vert t_2 - t_1 \Vert \ge \frac{1}{x}$, this already gives the claim. Otherwise, we can simply bound
    \[ \sum_{x < k \le \frac{1}{\Vert t_2 - t_1 \Vert}} \frac{\cos \big( k (t_2 - t_1) \big)}{2 k e^{2 k \sigma}} \le \frac{1}{2} \big( \log 1 / \Vert t_2 - t_1 \Vert - \log x \big) + O(1), \]
    which completes the proof.
\end{proof}

\section{Probabilistic tools}

We begin by recording the following results, which are simply Probability Results $1$ and $2$ from \cite{Harper1} and Probability Result $1$ from \cite{Harper6} as well as Normal Comparison Result $1$ from \cite{Harper8}.

\begin{proposition}[Harper]\label{PR1}
    Let $a \ge 1$. For any integer $n \ge 1$, let $G_1, \dots, G_n$ be independent real Gaussian random variables, each having mean zero and variance between $\frac{1}{20}$ and $20$, say. Let $h$ be a function such that $|h(j)| \le 10 \log j$. Then we have
    \[ \mathbb{P} \bigg[ \sum_{m = 1}^j G_m \le a + h(j) \, \forall \, 1 \le j \le n \bigg] \asymp \min \Big \{ 1, \frac{a}{\sqrt{n}} \Big \}. \]
\end{proposition}

\begin{proposition}[Harper]\label{PR2}
    There is an absolute constant $B$ such that the following is true. Let $a$ and $n$ be large, $n$ being an integer, and let $G_1, \dots G_n$ be Gaussian random variables, each with mean between $\frac{1}{20}$ and $20$, say. Then uniformly for any functions $h(j)$ resp. $g(j)$ satisfying $|h(j)| \le 10 \log j$ resp. $g(j) \le -B j$ for all $1 \le j \le n$, we have
    \[ \mathbb{P} \Big[ g(j) \le \sum_{m=1}^j G_m \le \min \{a, Bj \} + h(j) \, \forall \, 1 \le j \le n \Big] \asymp \min \Big \{ 1, \frac{a}{\sqrt{n}} \Big \}. \]
\end{proposition}

The reader should think of this barrier condition essentially as a barrier of the shape $\sum_{m=1}^j G_m \le a$ for all $1 \le j \le n$. The random walk typically varies on a scale of $\sqrt{j}$ , so that the imposed lower bound is essentially irrelevant; this more general shape is only necessary for technical reasons to be explained later. On the other hand, the shift by $h(j)$ is so small compared to the typical variation of the random walk that we should not expect it to make a big difference either, as is indeed reflected by the Proposition. The added flexibility will nonetheless be rather useful for our purposes. 

\begin{proposition}\label{PR3}
    Let $a, b \in \mathbb{R}$ and $n \in \mathbb{N}$ all be large, and let $G_1, \dots, G_n$ be independent real Gaussian random variables, each having mean zero and variance between $\frac{1}{20}$ and $20$, say. Then we have the uniform upper bound
    \[ \mathbb{P} \bigg[ \sum_{m = 1}^j G_m \le a \, \forall \, 1 \le j \le n, \quad a-b \le \sum_{m = 1}^n G_m \le a \bigg] \ll \min \Big \{ 1, \frac{a}{\sqrt{n}} \Big \} \min \Big \{ 1, \frac{b}{\sqrt{n}} \Big \}^2. \]
\end{proposition}

While the previous statements all concerned barrier conditions on random walks, the last result will pertain to estimating the maximum of a collection of Gaussian random variables with small correlations. It generalises the fairly well-known result that the maximum of $n$ independent standard Gaussians is close to $\sqrt{2 \log n}$ with high probability.

\begin{proposition}\label{PR4}
    Suppose that $n \ge 2$, and that $\varepsilon > 0$ is sufficiently small. Let $G_1, \dots, G_n$ be mean $0$, variance $1$, real, jointly Gaussian random variables, and suppose that $\mathbb{E}[G_i G_j] \le \varepsilon$ whenever $i \ne j$. Then for any $100 \varepsilon \le \delta \le \frac{1}{100}$ (say), we have
    \begin{equation}
        \mathbb{P} \bigg[ \max_{1 \le i \le n} G_i \le \sqrt{(2 - \delta) \log n} \bigg] \ll \exp \left( - \Theta \left( \frac{n^{\delta/20}}{\sqrt{\log n}} \right) \right) + n^{-\delta^2/50 \varepsilon}.
    \end{equation}
\end{proposition}

For our analysis we will require a lower bound for $\int_0^{2 \pi} \big| F_N( e^{-\sigma + i t} ) \big|^2 \, dt$ that holds with high probability. For this purpose, we have the following

\begin{proposition}\label{MCR4}
    Uniformly for all large $N$, all $0 \le \sigma \le \frac{1}{N^{1/100}}$ and all $1 \le W \le (\log N)^{1/3}$ (say), we have
    \[ \mathbb{P} \bigg[ \int_0^{2 \pi} \big| F_N( e^{-\sigma + i t} ) \big|^2 \, dt \ge \frac{1}{e^{21W/10}} \frac{\min \{ N, \frac{1}{\sigma} \} }{\sqrt{\log N}} \bigg] \ge 1 - O(e^{-W/10}). \]
\end{proposition}

In the end, the proof of this proposition will boil down to a somewhat more precise version of \cite[Proposition 11.1]{SoundZaman1}. Unlike there, it does not seem like we can require the fairly simple constant barrier condition that is imposed, but instead need a slightly modified upper bound and a very weak lower bound of the same type as in \cite[Multiplicative Chaos Result 4]{Harper8} in order to make the argument work.

Define $N_\sigma$ such that $\log N_\sigma$ is the largest integer $\le \log \min \{ N, \frac{1}{\sigma} \}$, and let $C$ be a (sufficiently large) constant such that $W+C$ is an integer. Let $M := \log N_\sigma - W - C$, and for $1 \le m \le M$ and $0 \le t \le 2 \pi$, set
\[ Z_t^{(W)}(m) = Z_{t,\sigma}^{(W)}(m) := Z_{t, \sigma}(W+C+m) = \sum_{e^{W+C+m-1} < k \le e^{W+C+m}} \frac{\Re X(k) e^{i t k}}{\sqrt{k} e^{k \sigma}} \] and $Z^{(W)}(m) := Z_0^{(W)}(m)$. Later on, we might shift by other quantities and keep the analogous notation. We also define $\sigma_{m,W}^2 := \sigma_{m+W+C}^2$ as well as $\rho_{m,t}^{(W)} := \rho_{m+W+C,t}$ and $\rho_{m}^{(W)} := \rho_{m+W+C}$. 

Further, let $\mathcal{L} \subseteq [0, 2 \pi]$ denote the (random) set of $t$ satisfying the barrier condition
\[ - C j \le  \sum_{m=1}^j \Big( Z_t^{(W)}(m) - 2 \sigma_{m,W}^2\Big) \le C - 3 \log(j) \]
for all $1 \le j \le M$. Also, set
\[ F^{(\mathrm{med})}(e^{-\sigma+it}) := \exp \bigg( \sum_{e^{W+C} < k \le N_\sigma} \frac{X(k)e^{i t k}}{\sqrt{k} e^{k\sigma}} \bigg), \]
so that \[|F^{(\mathrm{med})}(e^{-\sigma+it})| = \exp \Big( \sum_{m=1}^M Z_t^{(W)}(m) \Big).\]

Our basic strategy for proving Proposition \ref{MCR4} is to apply the following lower bounds
\begin{align}
        \int_0^{2 \pi} \big| F_N( e^{-\sigma + i t} ) \big|^2 \, dt &\ge \int_{\mathcal{L}} \big| F_N( e^{-\sigma + i t} ) \big|^2 \, dt \nonumber \\
        & \ge \frac{\bigg( \int_{\mathcal{L}} \big| F^{(\mathrm{med})}( e^{-\sigma + i t} ) \big|^2 \, dt \bigg)^2}{\int_{\mathcal{L}} | F^{(\mathrm{med})}( e^{-\sigma + i t} ) \big|^2 \exp \Big( - 2 \sum_{\substack{k \le e^{W+C} \text{ or} \\ N_\sigma < k \le N }} \frac{\Re X(k) e^{i t k}}{ \sqrt{k} e^{k \sigma} } \Big) \, dt  }, \label{eq31}
    \end{align}
where the second bound follows from Cauchy-Schwarz. Thus, we require an upper bound for the denominator here, as well as a lower bound for the numerator. The former will be a fairly direct matter to deal with, whereas the latter requires more careful analysis.

In order to obtain said lower bound for the numerator, the general strategy will be to apply Chebyshev's inequality in order to show that with high probability, it does not deviate too much from its mean (say at most half of the mean). Thus, we need an upper bound for the variance

\begin{align}
        &\mathbb{E} \Bigg[ \bigg( \int_{\mathcal{L}} \big| F^{(\mathrm{med})}( e^{-\sigma + i t} ) \big|^2 \, dt - \mathbb{E} \int_{\mathcal{L}} \big| F^{(\mathrm{med})}( e^{-\sigma + i t} ) \big|^2 \, dt \bigg)^2 \Bigg] \label{CorrelationBarrier} \\
        = 2 \pi \int_0^{2 \pi} &\mathbb{E} \Big[ \mathbbm{1}(0 \in \mathcal{L}) F^{(\mathrm{med})}( e^{-\sigma} ) \big|^2 \mathbbm{1}(t \in \mathcal{L}) F^{(\mathrm{med})}( e^{-\sigma + i t} ) \big|^2 \Big] \, dt - \mathbb{E} \bigg[ \int_{\mathcal{L}} \big| F^{(\mathrm{med})}( e^{-\sigma + i t} ) \big|^2 \, dt \bigg]^2. \nonumber
\end{align}

Hence, the next step is to acquire good bounds for the correlation
\begin{align}
    \mathbb{E} \Big[ \mathbbm{1}(0 \in \mathcal{L}) F^{(\mathrm{med})}( e^{-\sigma} ) \big|^2 \mathbbm{1}(t \in \mathcal{L}) F^{(\mathrm{med})}( e^{-\sigma + i t} ) \big|^2 \Big].
\end{align}
We will obtain various bounds depending on the distance $\Vert t \Vert$ of $t$ to the nearest multiple of $2 \pi$.

When $\Vert t \Vert$ is not too small, we show that the two quantities are almost uncorrelated.

\begin{proposition}\label{LargetCorrelation}
    Suppose that $\Vert t \Vert \ge e^{-W}$. Then we have
    \begin{align}
        &\mathbb{E} \Big[ \mathbbm{1}(0 \in \mathcal{L}) \big| F^{(\mathrm{med})}( e^{-\sigma} ) \big|^2 \mathbbm{1}(t \in \mathcal{L}) \big| F^{(\mathrm{med})}( e^{-\sigma + i t} ) \big|^2 \Big] \\ = \bigg( 1 + &O \bigg( \frac{1}{\Vert t \Vert e^{W}} \bigg) \bigg) \mathbb{E} \Big[ \mathbbm{1}(0 \in \mathcal{L}) \big| F^{(\mathrm{med})}( e^{-\sigma} ) \big|^2\Big] \mathbb{E} \Big[ \mathbbm{1}(t \in \mathcal{L}) \big| F^{(\mathrm{med})}( e^{-\sigma + i t} ) \big|^2 \Big]. \nonumber
    \end{align}
\end{proposition}

\begin{proof}
    First, we note that by Proposition \ref{RhoEstimate}, on this range of $t$ we have $\sigma_{m,W}^2 \asymp 1$
    as well as
    \[ |\rho_{m,t}^{(W)}| \le \frac{1}{10^4 \Vert t \Vert e^{W+m}} \le \frac{1}{10^4} \]
    for all $1 \le m \le M$, by taking $C$ sufficiently large.
    
    By independence of the $Z_t^{(W)}(m)$ for different $m$, plugging in the density of the bivariate Gaussian we have
    \begin{align*}
        &\mathbb{E} \Big[ \mathbbm{1}(0 \in \mathcal{L}) \big| F^{(\mathrm{med})}( e^{-\sigma} ) \big|^2 \mathbbm{1}(t \in \mathcal{L}) \big| F^{(\mathrm{med})}( e^{-\sigma + i t} ) \big|^2 \Big] \\
        &= \Big( \prod_{m=1}^M \frac{1}{2 \pi \sigma_{m,W}^2 \sqrt{1-\rho_m^{(W) 2}}} \Big) \idotsint \limits_{\mathcal{R}} \exp \bigg( \sum_{m=1}^M (2 x_m+ 2 y_m) - \sum_{m = 1}^M \frac{x_m^2 - 2 \rho_m^{(W)} x_m y_m + y_m^2}{2(1-\rho_m^{(W)2}) \sigma_{m,W}^2} \bigg) d \mathbf{x} \, d \mathbf{y},
    \end{align*}
    where
    \[ \mathcal{R} := \bigg \{ (\mathbf{x},\mathbf{y}) \in \mathbb{R}^M \times \mathbb{R}^M : - C j \le \sum_{m=1}^j \left( x_m - 2 \sigma_{m,W}^2 \right), \sum_{m=1}^j \left( y_m - 2 \sigma_{m,W}^2 \right) \le C - 3 \log(j) \; \forall \; 1 \le j \le M \bigg \}. \]
    Substituting $x_m' := x_m - 2 \sigma_{m,W}^2$ and $y_m' := y_m - 2 \sigma_{m,W}^2$, this equates to
    \begin{align*}
        &\exp \bigg( 4 \sum_{m = 1}^M \sigma_{m,W}^2 + 4 \sum_{m = 1}^M \frac{\rho_m^{(W)} \sigma_{m,W}^2}{1+\rho_m^{(W)}} \bigg) \Big( \prod_{m=1}^M \frac{1}{2 \pi \sigma_{m,W}^2 \sqrt{1-\rho_m^{(W)2}}} \Big) \times \\
        &\times \idotsint \limits_{\mathcal{R}'} \exp \bigg( - \sum_{m=1}^M \frac{x_m'^2 - 2 \rho_m^{(W)} x_m' y_m' + y_m'^2}{2 (1-\rho_m^{(W)2}) \sigma_{m,W}^2} + \sum_{m=1}^M \frac{2 \rho_m^{(W)} (x_m' + y_m')}{1 + \rho_m^{(W)}}  \bigg) d \mathbf{x}' \, d \mathbf{y}'
    \end{align*}
    with
    \[ \mathcal{R}' := \bigg \{ (\mathbf{x}',\mathbf{y}') \in \mathbb{R}^M \times \mathbb{R}^M : - C j \le \sum_{m=1}^j x_m', \sum_{m=1}^j y_m' \le C - 3 \log(j) \; \forall \; 1 \le j \le M \bigg \}. \]
    At this point we make crucial use of the weak lower bound that we implemented in the barrier condition, which implies in particular that $|x_m'|, |y_m'| \ll m$. Thus, this last expression is in turn
    \begin{align}
        &= \exp \bigg( O \bigg( \sum_{m=1}^M m^2 \rho_m^{(W)} \bigg) \bigg) \exp \bigg( 4 \sum_{m = 1}^M \sigma_{m,W}^2 \bigg) \bigg( \prod_{m=1}^M \frac{1}{2 \pi \sigma_{m,W}^2} \bigg) \idotsint \limits_{\mathcal{R}'} \exp \bigg( - \sum_{m=1}^M \frac{x_m'^2 + y_m'^2}{2 \sigma_m^2} \bigg) d \mathbf{x}' \, d \mathbf{y}' \nonumber \\
        &= \bigg( 1 + O \bigg( \frac{1}{\Vert t \Vert e^W} \bigg) \bigg) \exp \bigg( 2 \sum_{m = 1}^M \sigma_{m,W}^2 \bigg) \bigg( \prod_{m=1}^M \frac{1}{\sqrt{2 \pi \sigma_{m,W}^2}} \bigg) \idotsint \limits_{\mathcal{R}_1'} \exp \bigg( - \sum_{m=1}^M \frac{x_m'^2}{2 \sigma_{m,W}^2} \bigg) d \mathbf{x}' \times \nonumber \\
        &\qquad \qquad \qquad \qquad \quad \times \exp \bigg( 2 \sum_{m = 1}^M \sigma_{m,W}^2 \bigg) \bigg( \prod_{m=1}^M \frac{1}{\sqrt{2 \pi \sigma_{m,W}^2}} \bigg) \idotsint \limits_{\mathcal{R}_2'} \exp \bigg( - \sum_{m=1}^M \frac{y_m'^2}{2 \sigma_{m,W}^2} \bigg) d \mathbf{y}', \label{ProductOfExpectations}
    \end{align}
    where $\mathcal{R}_1'$ and $\mathcal{R}_2'$ describe the individual conditions for $x_m'$ resp. $y_m'$ from $\mathcal{R}'$. 
    
    On the other hand, we have
    \begin{align*}
    \mathbb{E} \Big[ \mathbbm{1}(0 \in \mathcal{L}) | F^{(\mathrm{med})}( e^{-\sigma} ) \big|^2 \Big] = \Big( \prod_{m=1}^M \frac{1}{\sqrt{2 \pi \sigma_{m,W}^2}} \Big) \idotsint \limits_{\mathcal{R}} \exp \left( \sum_{m=1}^M \Big( 2 z_m - \frac{z_m^2}{2 \sigma_{m,W}^2}  \Big) \right) \, d\mathbf{z},
    \end{align*}
    where \[ \mathcal{R} := \bigg \{ (z_m)_m :  - C j \le \sum_{m=1}^j \left( z_m - 2 \sigma_{m,W}^2 \right) \le C - 3 \log(j) \; \forall \; 1 \le j \le M \bigg \} \]
    (and the same at $t$ by translation invariance in law).
    Again substituting $z_m' = z_m - 2 \sigma_{m,W}^2$, this equates to
    \begin{align} \exp \bigg( 2 \sum_{m=1}^M  \sigma_{m,W}^2 \bigg) \prod_{m=1}^M \frac{1}{\sqrt{2 \pi \sigma_{m,W}^2}} \idotsint \limits_{\mathcal{R}'} \exp \bigg( - \sum_{m=1}^M \frac{z_m'^2}{2 \sigma_{m,W}^2} \bigg) \, d\mathbf{z}', \label{BarrierSecondMoment} \end{align}
    where
    \[ \mathcal{R}' := \bigg \{ (z_m')_m :  - C j \le \sum_{m=1}^j z_m' \le C - 3 \log(j) \; \forall \; 1 \le j \le M \bigg \}. \]
    Comparing \eqref{ProductOfExpectations} and \eqref{BarrierSecondMoment}, the claim follows. 
\end{proof}

We would like to record here for future use that by Proposition \ref{RhoEstimate} and \eqref{BarrierSecondMoment} we have
    \begin{align}
        \mathbb{E} \Big[ \mathbbm{1}(0 \in \mathcal{L}) | F^{(\mathrm{med})}( e^{-\sigma} ) \big|^2 \Big] \asymp \exp \Big( \sum_{m = 1}^M (1+O(e^{-W-m}) \Big) \mathbb{P}[\mathcal{B}] \asymp \frac{N_\sigma}{e^W} \mathbb{P}[\mathcal{B}],
    \end{align}
where $\mathcal{B}$ denotes the event that a sequence of Gaussians $G_1, \dots G_M$ with means $0$ and variances $\sigma_{m,W}^2$ satisfies the barrier condition defined by $\mathcal{R}'$, i.e. that for all $1 \le j \le M$ we have
\[ - C j \le \sum_{m = 1}^j G_m \le C - 3 \log j. \]
Thus, Proposition \ref{PR2} implies that
\begin{equation}
    \label{OneDimProbL}
    \mathbb{E} \Big[ \mathbbm{1}(0 \in \mathcal{L}) | F^{(\mathrm{med})}( e^{-\sigma} ) \big|^2 \Big] \asymp \frac{N_\sigma}{e^W \sqrt{\log N_\sigma}}.
\end{equation}

Next, we will deal with the range $0 < \Vert t \Vert \le e^{-W}$.

\begin{proposition}\label{SmalltCorrelation}
    For $0 < \Vert t \Vert \le N_\sigma^{-1/3}$, we have
    \begin{align}\label{Smallestt}
        &\mathbb{E} \Big[ \mathbbm{1}(0 \in \mathcal{L}) \big| F^{(\mathrm{med})}( e^{-\sigma} ) \big|^2 \mathbbm{1}(t \in \mathcal{L}) \big| F^{(\mathrm{med})}( e^{-\sigma + i t} ) \big|^2 \Big] \ll \frac{N_\sigma^2}{e^{3 W} (\log N_\sigma)^6} \min \Big \{ N_\sigma, \frac{1}{\Vert t \Vert} \Big \}.
    \end{align}
    For $N_\sigma^{-1/3} \le \Vert t \Vert \le e^{-W}$, we have
    \begin{align}\label{SecondSmallestt}
        &\mathbb{E} \Big[ \mathbbm{1}(0 \in \mathcal{L}) \big| F^{(\mathrm{med})}( e^{-\sigma} ) \big|^2 \mathbbm{1}(t \in \mathcal{L}) \big| F^{(\mathrm{med})}( e^{-\sigma + i t} ) \big|^2 \Big] \nonumber \\ 
        &\qquad \qquad \qquad \qquad \qquad \qquad \qquad \ll \frac{N_\sigma^2 }{e^{3 W} \Vert t \Vert \log N_\sigma \big( \log 1/\Vert t\Vert - W + 1 \big)^2} .
    \end{align}
    
\end{proposition}

\begin{proof}
    We begin with the proof of \eqref{Smallestt}, and make use of the parameter
    \[ M_1 :=\min \big \{ M, \big [\log (1/\Vert t \Vert) \big] - [W] + 1 \big \}, \]
    which satisfies $1 \le M_1 \le M$.
    The purpose of $M_1$ is to distinguish the values of $m$ where $Z_m^{(W)}(0)$ and $Z_m^{(W)}(t)$ are strongly correlated from those where they are almost uncorrelated. Our barrier condition for $t$ implies that
    \begin{equation}\label{BarrierSmallK} \exp \bigg( 2 \sum_{m=1}^{M_1} Z_m^{(W)}(t) \bigg) \ll \exp \bigg( 4 \sum_{m=1}^{M_1} \sigma_{m,W}^2 \bigg) \frac{1}{M_1^4} \ll \frac{e^{2 M_1}}{M_1^4}. \end{equation}
    
    Note that since $\Vert t \Vert \le N_\sigma^{-1/3}$ we have $M_1 \asymp \log N_\sigma$.
    We therefore obtain, discarding the indicator functions after applying \eqref{BarrierSmallK}, that
    \begin{align*}
        &\mathbb{E} \Big[ \mathbbm{1}(0 \in \mathcal{L}) \big| F^{(\mathrm{med})}( e^{-\sigma} ) \big|^2 \mathbbm{1}(t \in \mathcal{L}) \big| F^{(\mathrm{med})}( e^{-\sigma + i t} ) \big|^2 \Big] \\
        &\ll \frac{e^{2 M_1}}{M_1^6} \mathbb{E} \bigg[ \exp \bigg( 2 \sum_{m=1}^M Z^{(W)}(m) + 2 \sum_{m=M_1+1}^M Z_t^{(W)}(m) \bigg) \bigg] \\ &= \frac{e^{2 M_1}}{M_1^6} \mathbb{E} \bigg[ \exp \bigg( 2 \sum_{m=1}^{M_1} Z^{(W)}(m) \bigg) \bigg] \mathbb{E} \bigg[ \exp \bigg( 2 \sum_{m=M_1+1}^M Z^{(W)}(m) + 2 \sum_{m=M_1+1}^M Z_t^{(W)}(m) \bigg) \bigg].
    \end{align*}
    Using Proposition \ref{RhoEstimate} and \ref{MomentsZ}, this is simply
    \[ \ll \frac{e^{2 M_1}}{M_1^6} e^{M_1} e^{2(M-M_1)} = \frac{e^{M_1+2M}}{M_1^6} \ll \frac{1}{(\log N_\sigma)^6} \frac{\min \{ N_\sigma, \frac{1}{\Vert t \Vert} \}}{e^W} \frac{N_\sigma^2}{e^{2W}} = \frac{N_\sigma^2}{e^{3 W} (\log N_\sigma)^6} \min \Big\{ N_\sigma, \frac{1}{\Vert t \Vert} \Big\}. \]
    We would like to point out here that we have made crucial use of the small shift $- 3 \log j$ in our barrier condition, which has produced this additional factor $(\log N_\sigma)^6$ in the denominator. Without this, our bound on the contribution from small $\Vert t \Vert$ would turn out to be insufficient. This will be even more apparent on the next range.
    
    Now we shall deal with the range $N_\sigma^{-1/3} \le \Vert t \Vert \le e^{-W}$. This time, we will modify $M_1$ slightly by setting
    \[ M_2 := M_1 + 2 [\log M_1] = \big [\log (1/\Vert t \Vert) \big] - [W] + 1 + 2 \Big[ \log \Big( \big [\log (1/\Vert t \Vert) \big] - [W] + 1 \Big)  \Big]. \]
    This now allows us to distinguish between weakly and strongly correlated $Z_t^{(W)}$ variables with some margin.
    
    Here we have $1 \le M_2 \le \frac{1}{2} \log N_\sigma$ (say) and hence $M-M_2 \asymp \log N_\sigma$. We need to exploit more information from our barrier condition. Note that it implies in particular that for all $M_2+1 \le j \le M$ we have
    \[ -C(j+1) + 3 \log M_2 \le \sum_{m=M_2+1}^j (Z^{(W)}(m) - \sigma_{m,W}^2) \le C(M_2+1) - 3 \log j \]
    and the same event for $Z_t^{(W)}(m)$ in place of $Z^{(W)}(m)$. Denoting this event (at both $0$ and $t$) by $\mathcal{L}'$, we obtain
    \begin{align*}
        &\mathbb{E} \Big[ \mathbbm{1}(0 \in \mathcal{L}) \big| F^{(\mathrm{med})}( e^{-\sigma} ) \big|^2 \mathbbm{1}(t \in \mathcal{L}) \big| F^{(\mathrm{med})}( e^{-\sigma + i t} ) \big|^2 \Big] \\
        &\ll \frac{e^{3M_2}}{M_2^6} \mathbb{E} \bigg[ \mathbbm{1}(\mathcal{L}') \exp \bigg( 2 \sum_{m=M_2+1}^M Z^{(W)}(m) + 2 \sum_{m=M_2+1}^M Z_t^{(W)}(m) \bigg) \bigg].
    \end{align*}
    But in order to estimate this, we can employ a similar strategy as in the proof of Proposition \ref{LargetCorrelation}. Note that for any $1 \le m' \le M-M_2$, we have
    \[ |\rho_{m'}^{(W+M_2)}| \le \frac{1}{10^4 \Vert t \Vert e^{m'+W+M_2}} \le \frac{1}{10^3}, \]
    say. Plugging in the bivariate Gaussian density to obtain a $2(M-M_2)$-fold integral and making the same substitution, we arrive at
    \begin{align*}
        &\exp \bigg( 4 \sum_{m' = 1}^{M-M_2} \sigma_{m,W+M_2}^2 + 4 \sum_{m = 1}^{M-M_2} \frac{\rho_m^{(W+M_2)} \sigma_{m,W+M_2}^2}{1+\rho_m^{(W+M_2)}} \bigg) \Big( \prod_{m=1}^{M-M_2} \frac{1}{2 \pi \sigma_{m,W+M_2}^2 \sqrt{1-\rho_m^{(W+M_2)2}}} \Big) \times \\
        &\times \idotsint \limits_{\mathcal{R}'} \exp \bigg( - \sum_{m'=1}^{M-M_2} \frac{x_{m'}'^2 - 2 \rho_{m'}^{(W+M_2)} x_{m'}' y_{m'}' + y_{m'}'^2}{2 (1-\rho_{m'}^{(W+M_2)2}) \sigma_{m',W+M_2}^2} + \sum_{m'=1}^{M-M_2} \frac{2 \rho_{m'}^{(W+M_2)} (x_{m'}' + y_{m'}')}{1 + \rho_{m'}^{(W+M_2)}}  \bigg) d \mathbf{x}' \, d \mathbf{y}'
    \end{align*}
    with
    \begin{align*} &\mathcal{R}' := \bigg \{ (\mathbf{x}',\mathbf{y}') \in \mathbb{R}^{M-M_2} \times \mathbb{R}^{M-M_2} : \\ &- C (j+M_2+1) + 3 \log M_2 \le \sum_{m'=1}^j x_{m'}', \sum_{m'=1}^j y_{m'}' \le C(M_2+1) - 3 \log(j+M_2) \; \forall \; 1 \le j \le M - M_2 \bigg \}. \end{align*}
    Therefore, we have $|x_{m'}'|, |y_{m'}'| \ll m'+M_2$ and we can see that the previous expression is
    \begin{align*}
        &= \exp \bigg( O \bigg( \sum_{m'=1}^{M-M_2} (m'+M_2)^2 |\rho_{m'}^{(W+M_2)}| \bigg) \bigg) \exp \bigg( 4 \sum_{m = 1}^{M-M_2} \sigma_{m,W+M_2}^2 \bigg) \bigg( \prod_{m=1}^{M-M_2} \frac{1}{2 \pi \sigma_{m,W+M_2}^2} \bigg) \times \\  &\qquad \qquad \qquad \qquad \qquad \qquad \qquad \qquad \times \idotsint \limits_{\mathcal{R}'} \exp \bigg( - \sum_{m'=1}^{M-M_2} \frac{x_{m'}'^2 + y_{m'}'^2}{2 \sigma_{m',W+M_2}^2} \bigg) d \mathbf{x}' \, d \mathbf{y}' \\
        &\ll \exp \bigg( O \bigg( \frac{M_2^2}{\Vert t \Vert e^{W+M_2}} \bigg) \bigg) e^{2(M-M_2)} \mathbb{P}[\mathcal{B}]^2 \ll e^{2(M-M_2)} \mathbb{P}[\mathcal{B}]^2.
    \end{align*}
    Here, $\mathcal{B}$ is the event that a sequence of Gaussians $G_1, \dots, G_{M-M_2}$ with mean $0$ and variance $\sigma_{1,W+M_2}^2, \dots \sigma_{M-M_2,W+M_2}^2$ satisfies that for all $1 \le j \le M - M_2$ we have
    \[ \sum_{m'=1}^j G_{m'} \le C(M_2+1) - 3 \log(j+M_2). \]
    Moreover, we have used Proposition \ref{RhoEstimate} in the second step, and in the last step we crucially make use of the definition of $M_2$ (with our old parameter $M_1$, this error would not be acceptable).
    In order to understand the event $\mathcal{B}$, we can employ Proposition \ref{PR1} with $a = C(M_2+1) - 3 \log(M_2+1)$ and $n = M - M_2$ to deduce that $\mathbb{P}[\mathcal{B}] \ll \frac{M_2}{\sqrt{M-M_2}}$. We therefore conclude that
    \begin{align*}
        &\mathbb{E} \Big[ \mathbbm{1}(0 \in \mathcal{L}) \big| F^{(\mathrm{med})}( e^{-\sigma} ) \big|^2 \mathbbm{1}(t \in \mathcal{L}) \big| F^{(\mathrm{med})}( e^{-\sigma + i t} ) \big|^2 \Big] \\
        &\ll \frac{e^{3M_2}}{M_2^6} e^{2(M-M_2)} \frac{M_2^2}{M-M_2} = \frac{e^{2 M + M_2}}{(M-M_2) M_2^4} \ll \frac{N_\sigma^2 }{e^{3 W} \Vert t \Vert \log N_\sigma \big( \log 1/\Vert t\Vert - W + 1 \big)^2}
    \end{align*}
    as claimed.
\end{proof}

We are now in a position to complete the proof of Proposition \ref{MCR4}.

\begin{proof}[Conclusion of the proof of Proposition \ref{MCR4}.]
    Most of the work has already been done and we essentially just need to integrate our respective bounds over their ranges of $t$. 
    
    We begin by obtaining an upper bound for the denominator in \eqref{eq31} that holds with high probability. To this end, we note that by Proposition \ref{MomentsZ} and \eqref{OneDimProbL} its expectation satisfies
    \begin{align*}
        &\mathbb{E} \bigg[ \int_{\mathcal{L}} \big| F^{(\mathrm{med})}( e^{-\sigma + i t} ) \big|^2 \exp \Big( - 2 \sum_{\substack{k \le e^{W+C} \text{ or} \\ N_\sigma < k \le N }} \frac{\Re X(k) e^{i t k}}{ \sqrt{k} e^{k \sigma} } \Big) \, dt \bigg] \\
        &= \int_0^{2 \pi} \mathbb{E} \bigg[ \mathbbm{1}(t \in \mathcal{L}) \big| F^{(\mathrm{med})}( e^{-\sigma + i t} ) \big|^2 \bigg] \mathbb{E} \bigg[ \exp \Big( - 2 \sum_{\substack{k \le e^{W+C} \text{ or} \\ N_\sigma < k \le N }} \frac{\Re X(k) e^{i t k}}{ \sqrt{k} e^{k \sigma} } \Big) \bigg] \, dt \\
        &\ll e^W \mathbb{E} \Big[ \mathbbm{1}(0 \in \mathcal{L}) \big| F^{(\mathrm{med})}( e^{-\sigma} ) \big|^2 \Big] \ll \frac{N_\sigma}{\sqrt{\log N_\sigma}}.
    \end{align*}
    In the second step we have used that
    \[ \sum_{k > N_\sigma} \frac{1}{k e^{2 k \sigma}} \ll 1 \]
    to bound the contribution from the large $k$, which follows e.g. by dividing into dyadic intervals and bounding trivially on each of them.
    
    Hence, by Markov's inequality, with probability $\ge 1 - O(e^{-W/10})$ we have \begin{equation}\label{UpperBoundDenominator} \int_{\mathcal{L}} \big| F^{(\mathrm{med})}( e^{-\sigma + i t} ) \big|^2 \exp \Big( - 2 \sum_{\substack{k \le e^{W+C} \text{ or} \\ N_\sigma < k \le N }} \frac{\Re X(k) e^{i t k}}{ \sqrt{k} e^{k \sigma} } \Big) \, dt  \le \frac{e^{W/10} N_\sigma}{\sqrt{\log N_\sigma}}. \end{equation}
    
    As for the numerator, we note first that the contribution from small $\Vert t \Vert$ to the first integral in \eqref{CorrelationBarrier} is small: By Propositions \ref{LargetCorrelation} and \ref{SmalltCorrelation} as well as by \eqref{OneDimProbL} we have
    \begin{align*}
        &\int_{\Vert t \Vert \le e^{-W/2}} \mathbb{E} \Big[ \mathbbm{1}(0 \in \mathcal{L}) F^{(\mathrm{med})}( e^{-\sigma} ) \big|^2 \mathbbm{1}(t \in \mathcal{L}) F^{(\mathrm{med})}( e^{-\sigma + i t} ) \big|^2 \Big] \, dt \\
        &\ll \frac{N_\sigma^2}{e^{3 W} (\log N_\sigma)^6} \int_0^{N_\sigma^{-1/3}} \! \! \! \! \! \! \min \Big\{ N_\sigma, \frac{1}{t} \Big\} \, dt + \frac{N_\sigma^2}{e^{3 W} \log N_\sigma} \int_{N_\sigma^{-1/3}}^{e^{-W}} \frac{1}{t (\log 1/t - W + 1)^2} \, dt + \frac{N_\sigma^2}{e^{5W/2} \log N_\sigma} \\
        &\ll \frac{N_\sigma^2}{e^{3W} (\log N_\sigma)^5 } + \frac{N_\sigma^2}{e^{3 W} \log N_\sigma} + \frac{N_\sigma^2}{e^{5W/2} \log N_\sigma} \ll \frac{N_\sigma^2}{e^{5W/2}\log N_\sigma}.
    \end{align*}
    As a consequence we deduce, again invoking Proposition \ref{LargetCorrelation}, that the variance of the numerator in \eqref{eq31} (compare also \eqref{CorrelationBarrier}) satisfies
    \begin{align*}
        &\mathbb{E} \Bigg[ \bigg( \int_{\mathcal{L}} \big| F^{(\mathrm{med})}( e^{-\sigma + i t} ) \big|^2 \, dt - \mathbb{E} \int_{\mathcal{L}} \big| F^{(\mathrm{med})}( e^{-\sigma + i t} ) \big|^2 \, dt \bigg)^2 \Bigg] \\
        &\ll e^{-W/2} \mathbb{E} \Big[ \mathbbm{1}(0 \in \mathcal{L}) \big| F^{(\mathrm{med})}( e^{-\sigma} ) \big|^2 \Big]^2 + \frac{N_\sigma^2}{e^{5W/2}\log N_\sigma} \ll e^{-W/2} \mathbb{E} \Big[ \mathbbm{1}(0 \in \mathcal{L}) \big| F^{(\mathrm{med})}( e^{-\sigma} ) \big|^2 \Big]^2.
    \end{align*}
    Hence, by Chebyshev's inequality, we have that with probability $\ge 1 - O(e^{-W/2})$ the numerator in \eqref{eq31} satisfies
    \[ \int_{\mathcal{L}} \big| F^{(\mathrm{med})}( e^{-\sigma + i t} ) \big|^2 \, dt \ge \frac{1}{2} \mathbb{E} \bigg[ \int_{\mathcal{L}} \big| F^{(\mathrm{med})}( e^{-\sigma + i t} ) \big|^2 \, dt \bigg] \gg \frac{N_\sigma}{e^W \sqrt{\log N_\sigma}}. \]
    Combining this with \eqref{eq31} and \eqref{UpperBoundDenominator} (and the fact that $\log N_\sigma \asymp \log N$), we conclude that with probability $\ge 1 - O(e^{-W/10})$ we have
    \[ \int_0^{2 \pi} \big| F_N( e^{-\sigma + i t} ) \big|^2 \, dt \ge \frac{N_\sigma}{e^{21W/10} \sqrt{\log N}}, \]
    which gives the claim.
\end{proof}

\section{Reduction to a covariance estimate}

We continue to proceed in a rather similar fashion as in \cite{Harper8} and \cite{SoundZaman1}. The first step in the reduction is to discard the $N$-smooth part. For $N < n < 2 N$, we write
\begin{equation}\label{RoughPartOut} A(n) = \sum_{\lambda \pt n} a(\lambda) = \sum_{N < k \le n} \frac{X(k)}{\sqrt{k}} \sum_{\sigma \pt n - k} a(\sigma)  + \sum_{\substack{ \lambda \pt n \\ \lambda_1 \le N }} a(\lambda). \end{equation}
Here we used the fact that for $n$ in this range there can be at most one partition of length $ > N$, as well as that $a((k)) = \frac{X(k)}{\sqrt{k}}$, where $(k)$ denotes the partition of $k$ with only one part. As for the $N$-smooth part, note that
\begin{align*}
    &\mathbb{E} \bigg[ \, \bigg| \bigg\{ \frac{8 N}{7} \le n \le \frac{4N}{3} :\Big| \sum_{\substack{ \lambda \pt n \\ \lambda_1 \le N }} a(\lambda) \Big| \ge (\log N)^{1/100} \bigg\} \bigg| \, \bigg] \\ \le &\sum_{\frac{8 N}{7} \le n \le \frac{4N}{3}} \frac{\mathbb{E} \bigg[ \, \Big| \sum_{\substack{ \lambda \pt n \\ \lambda_1 \le N }} a(\lambda) \Big|^2 \bigg] }{(\log N)^{1/50}} \ll \frac{N}{(\log N)^{1/50}}
\end{align*}
using orthogonality and the fact that $\mathbb{E}[|A(n)|^2] = 1$ for all $n$. Thus, with probability $\ge 1 - O((\log N)^{-1/50}$ there exists a random set $\mathcal{N} \subseteq \Big[ \frac{8N}{7}, \frac{4N}{3} \Big]$ of cardinality  $\ge \frac{N}{7}$ (say) such that 
\begin{equation}\label{RoughPartBound} \Big| \sum_{\substack{ \lambda \pt n \\ \lambda_1 \le N }} a(\lambda) \Big| \le (\log N)^{1/100} \end{equation}
for all $n \in \mathcal{N}$. So far, our argument would have worked just as well on our original sum, but the crucial point is that the set $\mathcal{N}$ only depends on the random variables $(X(k))_{k \le N}$. Since we will now condition on these random variables, the set $\mathcal{N}$ will be fixed under this conditioning.

For $n \in \mathcal{N}$, let $Y_n$ be the random variables defined by
\[ \Re \sum_{N < k \le n} \frac{X(k)}{\sqrt{k}} \sum_{\sigma \pt n - k} a(\sigma)  \]
conditional on the values of $(X(k))_{k \le N}$,
so that $Y_n$ follows a Gaussian distribution, now being a linear combination of independent Gaussians. We will show that this real part is large enough with high probability, and from that deduce the claim for the absolute value.

We have
\[ \mathbb{E}[Y_n] = 0 \]
and
\begin{align*}
    \mathbb{E}[Y_n^2] &= \frac{1}{4} \mathbb{E} \bigg[ \bigg( \sum_{N < k \le n} \frac{X(k)}{\sqrt{k}} \sum_{\sigma \pt n - k} a(\sigma) + \sum_{N < k \le n} \frac{\overline{X(k)}}{\sqrt{k}} \sum_{\sigma \pt n - k} \overline{a(\sigma)} \bigg)^2 \, \Big \vert \, (X(k))_{k \le N} \bigg]  \\ &= \sum_{N < k \le n} \frac{1}{2k} \Big| \sum_{\sigma \pt n - k} a(\sigma) \Big|^2.
\end{align*}
In a similar manner, we obtain
\[ \mathbb{E} [Y_n Y_m] = \Re \sum_{N < k \le \frac{4N}{3}} \frac{1}{2k} \Big( \sum_{\sigma \pt n - k} a(\sigma) \Big) \Big( \sum_{\sigma \pt m - k} \overline{a(\sigma)} \Big), \]
where we interpret the respective sums as empty whenever $n-k$ resp. $m-k$ is negative.

The next step is to obtain lower bounds for the variance that hold with high probability.

\begin{proposition}\label{VarianceProposition}
Let $\frac{8N}{7} \le n \le \frac{4N}{3}$ and $1 \le W \le (\log N)^{1/3}$. Then with probability $\ge 1 - O(e^{-W/10})$, we have
\[ \sum_{N < k \le n} \frac{1}{k} \Big| \sum_{\sigma \pt n - k} a(\sigma) \Big|^2 \ge \frac{1}{e^{11W/5} \sqrt{\log N}}. \]
\end{proposition}

\begin{proof} We follow ideas from \cite[Proposition 1]{Harper8} and \cite[Proposition 8.1]{SoundZaman1}. 
Firstly, note that on the prescribed range of $k$ we have $\frac{1}{k} \ge \frac{3}{4N}$, and after replacing $k$ by $n-k$ we get 
\[ \sum_{N < k \le n} \frac{1}{k} \Big| \sum_{\sigma \pt n - k} a(\sigma) \Big|^2 \ge \frac{3}{4N} \sum_{k \le n-N} \Big| \sum_{\sigma \pt k} a(\sigma) \Big|^2 \gg \frac{1}{N} \sum_{k \le \frac{N}{7}} \Big| \sum_{\substack{ \sigma \pt k \\ \sigma_1 \le N}} a(\sigma) \Big|^2,  \]
where we inserted the void condition $\sigma_1 \le N$ for use in the following step.
Now for any $0 \le r \le 1$, we get that this expression is
\[ \ge \frac{1}{N} \sum_{k \ge 0} r^{7 k} \Big| \sum_{\substack{ \sigma \pt k \\ \sigma_1 \le N}} a(\sigma) \Big|^2 - \frac{r^N}{N} \sum_{k \ge 0} \Big| \sum_{\substack{ \sigma \pt k \\ \sigma_1 \le N}} a(\sigma) \Big|^2.  \]
Next, noting that
\[ F_N(z) = \sum_{k \ge 0} \bigg(\sum_{\substack{\sigma \pt k \\ \sigma_1 \le N}} a(\sigma) \bigg) z^k, \]
Parseval's identity implies that this equates to
\[ \frac{1}{2 \pi N} \int_0^{2 \pi} \left| F_N( e^{-7V/2N + i t}) \right|^2 \, dt - \frac{e^{-V}}{2 \pi N} \int_0^{2 \pi} \left| F_N(e^{i t}) \right|^2 \, dt \]
with $r = e^{-V/N}$ and $V \ge 0$. Applying \cite[Proposition 3.2]{SoundZaman1} with $q = \frac{2}{3}$ (and $r = 1$), we obtain that the subtracted second term satisfies
\[ \mathbb{E} \bigg[ \bigg( \frac{1}{2 \pi N} \int_0^{2 \pi} \left| F_N(e^{i t}) \right|^2 \, dt \bigg)^{2/3} \bigg] \ll \left(\frac{1}{\sqrt{\log N}}\right)^{2/3}.  \]
Hence Markov inequality tells us (for any $A \ge 1$) that with probability $\ge 1- O(A^{-2/3})$ we have
\[ \int_0^{2 \pi} \left| F_N(e^{i t}) \right|^2 \, dt \le \frac{A}{\sqrt{\log N}}. \]
As for the first term, Proposition \ref{MCR4} tells us that with probability $\ge 1 - O(e^{-W/10})$, 
\[ \int_0^{2 \pi} \left| F_N( e^{-7V/2N + i t}) \right|^2 \, dt \ge \frac{\min \{ 1, \frac{2}{7 V} \} }{e^{21W/10} \sqrt{\log N} }. \]
Choosing (say) $V = 6W$ and $A = e^{3W}$ and assuming without loss of generality that $W$ is large enough, the claim follows upon combining these estimates.

\end{proof}

The covariance estimate which we shall rely on for now and prove in the next section is the following

\begin{proposition}\label{CovarianceProposition}
    Let $N$ be sufficiently large, and for $n \in \big[ \frac{8N}{7}, \frac{4N}{3}\big]$ define
    \[ \mathcal{B}_n := \Big \{ m \in \Big[ \frac{8N}{7}, \frac{4N}{3} \Big] \; : \; \Big| \sum_{N < k \le \frac{4N}{3}} \frac{1}{k} \Big( \sum_{ \sigma \pt n-k} a(\sigma) \Big) \Big( \sum_{ \sigma \pt m-k} \overline{a(\sigma)} \Big) \Big| \ge \frac{1}{(\log N)^{4/5}} \Big \}. \]
    Then with probability $\ge 1 - O((\log N)^{-1/10})$, we have $\max_n |\mathcal{B}_n| \le (\log N)^{7/10}$.
\end{proposition}

Assuming this proposition, we are now in a position to prove Theorem \ref{T2}. Again, this proceeds similar to the proof of Theorem $2$ in \cite{Harper8}, though with minor simplifications.

\begin{proof}[Proof of Theorem \ref{T2} assuming Proposition \ref{CovarianceProposition}]
    Suppose without loss of generality that $W$ is sufficiently large. Note first that by \eqref{RoughPartOut} and \eqref{RoughPartBound} we have
    \begin{align*}
        &\mathbb{P} \bigg[ \max_{\frac{8 N}{7} \le n \le \frac{4N}{3}} |A(n)| \le \frac{(\log N)^{1/4}}{e^{1.2 W}} \bigg]
        \le \mathbb{P} \bigg[ \max_{\frac{8 N}{7} \le n \le \frac{4N}{3}} \Re A(n) \le \frac{2 (\log N)^{1/4}}{e^{6W/5}} - (\log N)^{1/100} \bigg] \\
        &\le \mathbb{P} \bigg[ \max_{n \in \mathcal{N}} \Re \sum_{N < k \le n} \frac{X(k)}{\sqrt{k}} \sum_{\sigma \pt n - k} a(\sigma) \le \frac{2 (\log N)^{1/4}}{e^{6W/5}} \bigg] + O((\log N)^{-1/50}).
    \end{align*}
    We can write this last probability as the expectation of the probability conditional on $(X(k))_{k \le N}$, which gives rise to the random variables $(Y_n)_{n \in \mathcal{N}}$. Defining $V := \min_{n \in \mathcal{N}} \mathbb{E}[Y_n^2]$ (which is really a conditional variance), Proposition \ref{VarianceProposition} implies that with probability $\ge 1 - O(e^{-W/10})$ we have
    \[ V \gg \frac{1}{e^{11W/5} \sqrt{\log N}}, \]
    so that under this event we have
    \begin{align*} \mathbb{P} \bigg[ \max_{n \in \mathcal{N}} Y_n \le \frac{2 (\log N)^{1/4}}{e^{6W/5}} \bigg] &\le \mathbb{P} \bigg[ \max_{n \in \mathcal{N}} \frac{Y_n}{\sqrt{\mathbb{E}[Y_n^2]}} \le \frac{2 (\log N)^{1/4}}{\sqrt{V} e^{6W/5}} \bigg]
    \le \mathbb{P} \bigg[ \max_{n \in \mathcal{N}} \frac{Y_n}{\sqrt{\mathbb{E}[Y_n^2]}} \le \sqrt{\frac{1}{2} \log N} \bigg], \end{align*}
    say, since $W$ is sufficiently large. Note that the random variables $Y_n/\sqrt{\mathbb{E}[Y_n^2]}$ are jointly Gaussian and each of them is normalized to mean $0$ and variance $1$. 
    
    Next, Proposition \ref{CovarianceProposition} implies that, under this same event on $V$, with probability $\ge 1 - O((\log N)^{-1/10})$ the following holds. There exists a subset $\mathcal{N}'$ of $\mathcal{N}$ of cardinality $\gg N^{3/10}$, such that for any $n \ne m \in \mathcal{N}'$ we have the correlation estimate
    \[ \mathbb{E} \Big[ \frac{Y_n Y_m}{\sqrt{\mathbb{E}[Y_n^2]} \sqrt{\mathbb{E}[Y_m^2]}} \Big] \le \frac{(\log N)^{-4/5}}{2 V} \ll \frac{e^{11W/5}}{(\log N)^{3/10}}. \]
    By our assumption on $W$, we can thus make this correlation small, say smaller than $\varepsilon$ for some sufficiently small constant $0 < \varepsilon < 10^{-4}$; this puts us into a position to apply Proposition \ref{PR4} with (say) $\delta = \sqrt{\varepsilon}$, which implies that
    \begin{align*}
        \mathbb{P} \bigg[ \max_{n \in \mathcal{N}} \frac{Y_n}{\sqrt{\mathbb{E}[Y_n^2]}} \le \sqrt{\frac{1}{2} \log N} \bigg] &\le \mathbb{P} \bigg[ \max_{n \in \mathcal{N}} \frac{Y_n}{\sqrt{\mathbb{E}[Y_n^2]}} \le \sqrt{(2 - \sqrt{\varepsilon}) \log ( \#\mathcal{N}' )} \bigg] \\
        &\ll \exp \left( - N^{\sqrt{\varepsilon}/100} \right) + N^{-3/500} \ll N^{-3/500}.
    \end{align*}
    Putting all of this together implies that we have
    \[\max_{\frac{8 N}{7} \le n \le \frac{4N}{3}} |A(n)| \le \frac{(\log N)^{1/4}}{e^{1.2 W}} \]
    with probability
    \[ \ge 1 - O((\log N)^{-1/50}) - O(e^{-W/10}) - O((\log N)^{-1/10}) - O(N^{-3/500}) = 1 - O(e^{-W/10})  \]
\end{proof}

\section{Proof of the covariance estimate}

We begin by noting that as a consequence of \eqref{Cauchy} (which holds even for negative $N$ if we set $A(N)$ to be $0$), we have that the covariance sum we are interested in satisfies, for any $K \ge \frac{N}{3}$, the relation
\begin{align*}
    &\sum_{N < k \le \frac{4N}{3}} \frac{1}{k} \Big( \sum_{ \sigma \pt n-k} a(\sigma) \Big) \Big( \sum_{ \sigma \pt m-k} \overline{a(\sigma)} \Big) \\
    &= \sum_{N < k \le \frac{4N}{3}} \frac{1}{(2 \pi i)^2 k} \int_{|z_1|=1} \int_{|z_2|=1} F_K(z_1) \overline{F_K(z_2)} \frac{z_1^{m-k+1}}{z_2^{n-k+1}} \, dz_1 \, dz_2 \\
    &= \frac{1}{4 \pi^2} \int_0^{2 \pi} \int_0^{2 \pi} F_K(e^{i t_1}) \overline{F_K(e^{i t_2})} e^{i t_1 (m+1) - i t_2 (n+1)} \sum_{N < k \le \frac{4N}{3}} \frac{e^{i(t_2-t_1) k}}{k} \, dt_1 \, dt_2.
\end{align*}
We shall take $\frac{N}{3} \le K \le N$ (say) such that $\log K$ is an integer. 

Our first task is to restrict, with high probability, to an integration range where $\Vert t_2 - t_1 \Vert$ is fairly small.

\begin{proposition}
    Let $\frac{N}{3} \le K \le N$ (say) be such that $K$ is an integer, and let $\frac{8 N}{7} \le n, m \le \frac{4 N}{3}$. With probability $\ge 1 - O((\log N)^{-25})$ we have
    \begin{align}
        &\sum_{N < k \le \frac{4N}{3}} \frac{1}{k} \Big( \sum_{ \sigma \pt n-k} a(\sigma) \Big) \Big( \sum_{ \sigma \pt m-k} \overline{a(\sigma)} \Big) \\
        &= \frac{1}{4 \pi^2} \! \! \! \! \! \! \! \! \iint \limits_{\substack{ [0, 2\pi] \times [0, 2 \pi] \\ \Vert t_2 - t_1 \Vert \le \frac{(\log N)^{100}}{N} }} \! \! \! \! \! \! \! \! \! \! \! F_K(e^{i t_1}) \overline{F_K(e^{i t_2})} e^{i t_1 (m+1) - i t_2 (n+1)} \! \! \! \! \sum_{N < k \le \frac{4N}{3}} \frac{e^{i(t_2-t_1) k}}{k} \, dt_1 \, dt_2 + O((\log N)^{-25}). \nonumber
    \end{align}
\end{proposition}

\begin{proof}
    We proceed by bounding the expectation on the range where $\Vert t_2 - t_1 \Vert > \frac{(\log N)^{100}}{N}$ and then applying Markov's inequality. Note that by Proposition \ref{NTR2} we have
    \[ \sum_{N < k \le \frac{4N}{3}} \frac{e^{i(t_2-t_1) k}}{k} \ll \frac{1}{\Vert t_2 - t_1 \Vert N}. \]
    Thus we can bound the expectation of the integral on this complementary range by
    \begin{equation}\label{NonDiagonalContr}
        \ll \frac{1}{N} \iint \limits_{\substack{ [0, 2\pi] \times [0, 2 \pi] \\ \Vert t_2 - t_1 \Vert > \frac{(\log N)^{100}}{N} }} \! \! \! \! \! \! \! \! \! \frac{ \mathbb{E} \big[ |F_K(e^{i t_1})| |F_K(e^{i t_2})|\big]}{\Vert t_2 - t_1 \Vert} \, dt_1 \, dt_2.
    \end{equation}
    Then by Propositions \ref{MomentsZ} and \ref{NTR2} the expectation in the numerator satisfies
    \begin{align*}
        \mathbb{E} \big[ |F_K(e^{i t_1})| |F_K(e^{i t_2})|\big] &= 
        \exp \Big( \sum_{k \le K} \frac{1 + \cos(k(t_2-t_1))}{2 k} \Big)
        \ll \sqrt{N} \exp \Big( \Re \sum_{k \le K} \frac{e^{i (t_2 - t_1) k}}{2 k} \Big) \\
        &\ll \frac{\sqrt{N}}{\sqrt{\Vert t_2 - t_1 \Vert}} \exp \Big( \Re \sum_{ \frac{1}{\Vert t_2 - t_1 \Vert} < k \le K} \frac{e^{i (t_2 - t_1) k}}{2 k} \Big) \ll \frac{\sqrt{N}}{\sqrt{\Vert t_2 - t_1 \Vert}}.
    \end{align*}
    We can therefore bound \eqref{NonDiagonalContr} by
    \[ \ll \frac{1}{\sqrt{N}} \int_0^{2 \pi} \int_{\frac{(\log N)^{100}}{N}}^{\pi} \frac{1}{t^{3/2}} \, dt \, dv \ll \frac{1}{(\log N)^{50}}. \]
    Hence, by Markov's inequality, with probability $\ge 1 - O((\log N)^{-25})$ we have
    \[ \frac{1}{4 \pi^2} \! \! \! \! \! \! \! \! \iint \limits_{\substack{ [0, 2\pi] \times [0, 2 \pi] \\ \Vert t_2 - t_1 \Vert > \frac{(\log N)^{100}}{N} }} \! \! \! \! \! \! \! \! \! \! \! F_K(e^{i t_1}) \overline{F_K(e^{i t_2})} e^{i t_1 (m+1) - i t_2 (n+1)} \! \! \! \! \sum_{N < k \le \frac{4N}{3}} \frac{e^{i(t_2-t_1) k}}{k} \, dt_1 \, dt_2 = O((\log N)^{-25}) \]
    as claimed.
\end{proof}

The next step is to restrict to the high-probability event that a certain barrier condition holds. Namely, denote 
\[ \mathcal{D}_t := \Big \{ \sum_{m = 1}^j (Z_t(m) - 2 \sigma_m^2) \le 12 \log \log N \, \; \; \forall \, 1 \le j \le \log K \Big \}. \]
Then it follows directly from \cite[Proposition 5.2]{SoundZaman1}, taking $A = 2 \log \log N \le \sqrt{\log K}$ and noting $10 \log j \le 10 \log \log N$ for all $j$, that
\[ \mathbb{P}[\mathcal{D}_t \text{ fails for any } t] \ll e^{-(\log \log N)^2}. \]
We can therefore deduce that with probability $\ge 1 - O((\log N)^{-25})$ we in fact have

\begin{align}
        &\sum_{N < k \le \frac{4N}{3}} \frac{1}{k} \Big( \sum_{ \sigma \pt n-k} a(\sigma) \Big) \Big( \sum_{ \sigma \pt m-k} \overline{a(\sigma)} \Big) \label{CovarianceRestrictD} \\
        &= \frac{1}{4 \pi^2} \! \! \! \! \! \! \! \! \iint \limits_{\substack{ [0, 2\pi] \times [0, 2 \pi] \\ \Vert t_2 - t_1 \Vert \le \frac{(\log N)^{100}}{N} }} \! \! \! \! \! \! \! \! \! \! \! \mathbbm{1}(\mathcal{D}_{t_1}) F_K(e^{i t_1}) \mathbbm{1}(\mathcal{D}_{t_2}) \overline{F_K(e^{i t_2})} e^{i t_1 (m+1) - i t_2 (n+1)} \! \! \! \! \sum_{N < k \le \frac{4N}{3}} \frac{e^{i(t_2-t_1) k}}{k} \, dt_1 \, dt_2 \nonumber \\ &\qquad \qquad \qquad \qquad \qquad \qquad \qquad \qquad \qquad \qquad \qquad \qquad \qquad \qquad + O((\log N)^{-25}). \nonumber
    \end{align}
    
    Next, we want to restrict further to a strengthened barrier condition that may not hold with high probability, but the contribution to the above integral under the complementary event is nonetheless very small. Namely, set
    \[ \mathcal{A}_t := \Big \{ \sum_{m = 1}^j (Z_t(m) - 2 \sigma_m^2) \le - 2000 \log \log N \, \; \; \forall \, 0.01 \log K \le j \le \log K \Big \} \cap \mathcal{D}_t. \]
    In order to restrict to this event, we will prove the following analogue of \cite[Multiplicative Chaos Result 3]{Harper8}. 
    
    By a slight abuse of notation, from now on we will write $\sigma_{m,J}^2 := \sigma_{m+J}^2$, i.e. we are not shifting by an additional sufficiently large constant anymore.
    
    \begin{proposition}
        For any $0 \le t \le 2 \pi$ we have
        \begin{equation}
            \mathbb{E} \big[ \mathbbm{1}(\mathcal{D}_t) \mathbbm{1}(\mathcal{A}_t \text{ fails}) |F_K(e^{i t})|^2 \big] \ll \frac{N (\log \log N)^4}{\log N}.
        \end{equation}
    \end{proposition}
    
    \begin{proof}
        The proof proceeds along similar, though somewhat simpler lines since the discretisation step is hidden in the proof of \cite[Proposition 5.2]{SoundZaman1} and, as usual, we do not need any Gaussian approximation statements or prime number estimates.
        
        Note that if $\mathcal{D}_t$ holds but $\mathcal{A}_t$ fails then there must exist some $0.01 \log K \le J \le \log K$ such that
        \[ \sum_{m = 1}^J (Z_t(m) - 2 \sigma_m^2) > - 2000 \log \log N. \] We can therefore bound
        \[ \mathbb{E} \big[ \mathbbm{1}(\mathcal{D}_t) \mathbbm{1}(\mathcal{A}_t \text{ fails}) |F_K(e^{i t})|^2 \big] \le \sum_{0.01 \log K \le J \le \log K} \mathbb{E} \big[ \mathbbm{1}(\mathcal{A}_t(J)) \mathbbm{1}(\mathcal{B}_t(J)) |F_K(e^{i t})|^2 \big], \]
        where
        \[ \mathcal{A}_t(J) := \Big \{ \sum_{m = 1}^j (Z_t(m) - 2 \sigma_m^2) \le 12 \log \log N \; \; \forall 1 \le j \le J, \quad \sum_{m = 1}^J (Z_t(m) - 2 \sigma_m^2) > - 2000 \log \log N \Big \} \]
        and
        \[ \mathcal{B}_t(J) := \Big \{ \sum_{m = J+1}^j (Z_t(m) - 2 \sigma_m^2) \le 2012 \log \log N \; \; \forall J+1 \le j \le \log K \Big \}. \]
        Since $\mathcal{A}_t(J)$ only concerns the random variables $Z_t(m)$ with $m \le J$ and $\mathcal{B}_t(J)$ only concerns those with $m > J$, we can factor
        \[ \mathbb{E} \big[ \mathbbm{1}(\mathcal{A}_t(J)) \mathbbm{1}(\mathcal{B}_t(J)) |F_K(e^{i t})|^2 \big] = \mathbb{E} \big[ \mathbbm{1}(\mathcal{A}_t(J)) |F_{e^J}(e^{i t})|^2 \big] \mathbb{E} \Big[ \mathbbm{1}(\mathcal{B}_t(J)) \Big| \frac{F_K}{F_{e^J}}(e^{i t}) \Big|^2 \Big]. \]
        For the second expectation, we can invoke the same idea as in the proof of Proposition \ref{LargetCorrelation} of writing the expectation as a $(\log K - J)$-fold integral over variables $x_{1}, \dots, x_{\log K - J}$ and substituting $x_m' := x_m - 2 \sigma_{m,J}^2$ and we obtain that
        \[ \Big[ \mathbbm{1}(\mathcal{B}_t(J)) \Big| \frac{F_K}{F_{e^J}}(e^{i t}) \Big|^2 \Big] = \exp \Big( 2 \sum_{m = 1}^{\log K-J} \sigma_{m,J}^2) \Big) \mathbb{P}[B(J)], \]
        where $B(J)$ is the event that a sequence of Gaussians $(G_m)_{m=1}^{\log K - J}$ with mean $0$ and variances $\sigma_{m,J}^2$ (all between $\frac{1}{20}$ and $20$) satisfy
        \[ \sum_{m = 1}^j G_m \le 2012 \log \log N \]
        for all $1 \le j \le \log K - J$. By Proposition \ref{PR1}, we know that $\mathbb{P}[B(J)] \ll \frac{\log \log N}{\sqrt{\log K - J}+1}$ (this bound is trivial for $J = \log K)$. We therefore see, using Proposition \ref{MomentsZ}, that
        \[ \mathbb{E} \Big[ \mathbbm{1}(\mathcal{B}_t(J)) \Big| \frac{F_K}{F_{e^J}}(e^{i t}) \Big|^2 \Big] \ll \frac{K \log \log N}{e^J (\sqrt{\log K - J}+1)}. \]
        As for $\mathcal{A}_t(J)$, the same idea of writing the expectation as a $J$-fold integral and making the same substitution implies that
        \[ \mathbb{E} \big[ \mathbbm{1}(\mathcal{A}_t(J)) |F_{e^J}(e^{i t})|^2 \big] = \exp \Big( 2 \sum_{m = 1}^J \sigma_m^2 \Big) \mathbb{P}[A(J)], \]
        where $A(J)$ is the event that a sequence of Gaussians $(G_m)_{m=1}^{J}$ with mean $0$ and variances $\sigma_{m,J}^2$ (all between $\frac{1}{20}$ and $20$) satisfy
        \[ \sum_{m = 1}^j G_m \le 12 \log \log N \; \; \forall \, 1 \le j \le J, \qquad \sum_{m = 1}^J G_m > - 2000 \log \log N. \]
        So Proposition \ref{PR3} tells us that $\mathbb{P}[B(J)] \ll \frac{(\log \log N)^3}{J^{3/2}}$ and hence
        \[ \mathbb{E} \big[ \mathbbm{1}(\mathcal{A}_t(J)) |F_{e^J}(e^{i t})|^2 \big] \ll \frac{e^J (\log \log N)^3}{(\log K)^{3/2}}. \]
        Putting things together, we have shown that
        \begin{align*}
            \mathbb{E} \big[ \mathbbm{1}(\mathcal{D}_t) \mathbbm{1}(\mathcal{A}_t \text{ fails}) |F_K(e^{i t})|^2 \big] &\ll \sum_{0.01 \log K \le J \le \log K} \frac{K \log \log N}{e^J (\sqrt{\log K - J}+1)} \frac{e^J (\log \log N)^3}{(\log K)^{3/2}} \\ &\ll \frac{N (\log \log N)^4}{\log N}
        \end{align*}
        as claimed.
    \end{proof}
    
    Our next goal is to boost the estimate \eqref{CovarianceRestrictD} to an estimate under the stronger events $\mathcal{A}_{t_1}$ and $\mathcal{A}_{t_2}$ in place of $\mathcal{D}_{t_1}$ and $\mathcal{D}_{t_2}$. To this end, note that by Proposition \ref{NTR2} the error we pick up in doing so is
    \[ \ll \! \! \! \! \! \! \! \! \iint \limits_{\substack{ [0, 2\pi] \times [0, 2 \pi] \\ \Vert t_2 - t_1 \Vert \le \frac{(\log N)^{100}}{N} }} \! \! \! \! \! \! \! \! \! \! \! \mathbbm{1}(\mathcal{D}_{t_1}) \mathbbm{1}(\mathcal{A}_{t_1} \text{ fails}) |F_K(e^{i t_1})| \mathbbm{1}(\mathcal{D}_{t_2}) \mathbbm{1}(\mathcal{A}_{t_2} \text{ fails}) |F_K(e^{i t_2})| \min \Big \{ 1, \frac{1}{N \Vert t_2 - t_1 \Vert} \Big \} dt_1 \, dt_2. \]
    Using $|x y| \le x^2 + y^2$ (say) on $x = \mathbbm{1}(\mathcal{D}_{t_1}) \mathbbm{1}(\mathcal{A}_{t_1} \text{ fails}) |F_K(e^{i t_1})| $ and $y = \mathbbm{1}(\mathcal{D}_{t_2}) \mathbbm{1}(\mathcal{A}_{t_2} \text{ fails}) |F_K(e^{i t_2})|$, noting the symmetry of the two arising expressions, this term has expectation
    \[ \ll \int_0^{\frac{(\log N)^{100}}{N}} \mathbb{E} \big[ \mathbbm{1}(\mathcal{D}_t) \mathbbm{1}(\mathcal{A}_t \text{ fails}) |F_K(e^{i t})|^2 \big] \min \Big \{ 1, \frac{1}{N t} \Big \} \, dt \ll \frac{(\log \log N)^5}{\log N}. \]
    We therefore deduce that with probability $\ge 1 - O((\log N)^{-1/10})$, the error is $\le \frac{(\log \log N)^5}{(\log N)^{9/10}}$. Hence with probability $\ge 1 - O((\log N)^{-1/10})$, we have (say)
    \begin{align}
    &\sum_{N < k \le \frac{4N}{3}} \frac{1}{k} \Big( \sum_{ \sigma \pt n-k} a(\sigma) \Big) \Big( \sum_{ \sigma \pt m-k} \overline{a(\sigma)} \Big) \label{CovarianceRestrictA} \\
        &= \frac{1}{4 \pi^2} \! \! \! \! \! \! \! \! \iint \limits_{\substack{ [0, 2\pi] \times [0, 2 \pi] \\ \Vert t_2 - t_1 \Vert \le \frac{(\log N)^{100}}{N} }} \! \! \! \! \! \! \! \! \! \! \! \mathbbm{1}(\mathcal{A}_{t_1}) F_K(e^{i t_1}) \mathbbm{1}(\mathcal{A}_{t_2}) \overline{F_K(e^{i t_2})} e^{i t_1 (m+1) - i t_2 (n+1)} \! \! \! \! \sum_{N < k \le \frac{4N}{3}} \frac{e^{i(t_2-t_1) k}}{k} \, dt_1 \, dt_2 \nonumber \\ &\qquad \qquad \qquad \qquad \qquad \qquad \qquad \qquad \qquad \qquad \qquad \qquad \qquad \qquad + O \left( (\log N)^{-17/20} \right). \nonumber
        \end{align}
        
        \begin{proposition}\label{TediousProposition}
            For any large $N$ and any positive integer $2 \le p \le N^{1/4}$ (say), with probability $\ge 1 - O((\log N)^{-2 p})$ have
            \begin{align*}
                &\max_{\frac{8 N}{7} < m \le \frac{4N}{3}} \sum_{\frac{8 N}{7} < n \le \frac{4N}{3}} \bigg| \iint \limits_{\substack{ [0, 2\pi] \times [0, 2 \pi] \\ \Vert t_2 - t_1 \Vert \le \frac{(\log N)^{100}}{N} }} \! \! \! \! \! \! \! \! \! \! \! \mathbbm{1}(\mathcal{A}_{t_1}) F_K(e^{i t_1}) \mathbbm{1}(\mathcal{A}_{t_2}) \overline{F_K(e^{i t_2})} e^{i t_1 (m+1) - i t_2 (n+1)} \! \! \! \! \sum_{N < k \le \frac{4N}{3}} \frac{e^{i(t_2-t_1) k}}{k} \, dt_1 \, dt_2 \bigg|^{2 p} \\
                &\ll N^{2/3} \left( \frac{C \log \log N}{N} \int_0^{2 \pi} |F_K(e^{i t})|^2 \, dt \right)^{2 p} + N^{2/3} \left( p^{25} (\log N)^{27} \right)^{2 p} + N \left( \frac{p^{25}}{(\log N)^{985}} \right)^{2 p}.
            \end{align*}
        \end{proposition}
        
        Note that the restriction to $p \ge 2$ is not really necessary, but since we do not need the case $p = 1$ and since some of the indexing later in the proof would require a case distinction, we leave the proof of said case to the interested reader.
        
        \begin{proof}
            Aside from details, the proof of this Proposition is quite similar to that of \cite[Proposition 3]{Harper8}.
            
            The first step in the proof is to expand the $2 p$-th power to obtain a $2 p$-fold double integral, say over variables $u_1, \dots, u_{2 p}$ (corresponding to $t_1$) and $v_1, \dots, v_{2 p}$ (corresponding to $t_2$). While the arising expression is somewhat complicated, it should be noted that the only part in the above that depends on $n$ is the factor $e^{-i t_2(n+1)}$, which after expanding the $2 p$-th power gives rise to a factor $e^{-i (v_1 + \dots + v_p - v_{p+1} - \dots - v_{2 p}) (n+1)}$. We can then perform the sum over $n$, noting that
            \[ \sum_{\frac{8 N}{7} < n \le \frac{4 N}{3}} e^{-i (v_1 + \dots +  v_p - v_{p+1} - \dots - v_{2 p}) (n+1)} \ll \min \Big \{ N, \frac{1}{\Vert v_1 + \dots +  v_p - v_{p+1} - \dots - v_{2 p} \Vert } \Big \}.  \]
            One then takes absolute values inside the integral. Since the only dependence on $m$ in the original integral is in the factor $e^{i t_1 (m+1)}$, the arising expression in the expanded integral will have absolute value independent of $m$, so that we can dispose of the maximum. We arrive at the bound
            \begin{align*}
                &\max_{\frac{8 N}{7} < m \le \frac{4N}{3}} \sum_{\frac{8 N}{7} < n \le \frac{4N}{3}} \bigg| \iint \limits_{\substack{ [0, 2\pi] \times [0, 2 \pi] \\ \Vert t_2 - t_1 \Vert \le \frac{(\log N)^{100}}{N} }} \! \! \! \! \! \! \! \! \! \! \! \mathbbm{1}(\mathcal{A}_{t_1}) F_K(e^{i t_1}) \mathbbm{1}(\mathcal{A}_{t_2}) \overline{F_K(e^{i t_2})} e^{i t_1 (m+1) - i t_2 (n+1)} \! \! \! \! \sum_{N < k \le \frac{4N}{3}} \frac{e^{i(t_2-t_1) k}}{k} \, dt_1 \, dt_2 \bigg|^{2 p} \\
                &\ll \idotsint \limits_{u \in [0, 2\pi]^{2p}} \idotsint \limits_{\substack{v \in [0, 2\pi]^{2p}  \\ \Vert v_i - u_i \Vert \le \frac{(\log N)^{100}}{N}}} \prod_{i=1}^{2p} \bigg( \mathbbm{1}(\mathcal{A}_{u_i}) |F_K(e^{i u_i})| \mathbbm{1}(\mathcal{A}_{v_i}) |F_K(e^{i v_i})| \,  \bigg| \sum_{N < k \le \frac{4N}{3}} \frac{e^{i(v_i-u_i) k}}{k} \bigg| \bigg) \times \\
                &\qquad \qquad \qquad \qquad \times \min \Big \{ N, \frac{1}{\Vert v_p + \dots +  v_p - v_{p+1} - \dots - v_{2 p} \Vert } \Big \} \, dv_1 \dots dv_{2p} \, du_1 \dots du_{2 p}.
            \end{align*}
            
            In order to deal with this integral, we begin by analysing the part where $\Vert v_1 + \dots +  v_p - v_{p+1} - \dots - v_{2 p} \Vert \ge N^{-2/3}$, say. On this range, we shall bound the minimum simply by $N^{2/3}$, and we also discard the indicators. Moreover, we bound
            \[ \prod_{i=1}^{2p} |F_K(e^{i u_i})| |F_K(e^{i v_i})| \le \prod_{i=1}^{2p} |F_K(e^{i u_i})|^2 + \prod_{i=1}^{2p} |F_K(e^{i v_i})|^2. \]
            These bounds imply that this portion of the integral is
            \[ \ll N^{2/3} \bigg( \iint \limits_{\substack{ [0, 2\pi] \times [0, 2 \pi] \\ \Vert t_2 - t_1 \Vert \le \frac{(\log N)^{100}}{N} }} \! \! \! \! \! \! \! \! \! |F_K(e^{i t_1})|^2 \, \bigg| \sum_{N < k \le \frac{4N}{3}} \frac{e^{i(t_2-t_1) k}}{k} \bigg| \, dt_1 \, dt_2 \bigg)^{2 p}. \]
            Further, applying Proposition \ref{NTR2} to the inner sum (on the range $\Vert t_2 - t_1 \Vert \ge \frac{1}{N}$, otherwise the trivial bound) implies that this is in turn
            \[ \ll N^{2/3} \bigg( \frac{C \log \log N}{N} \int_0^{2 \pi} | F_K(e^{i t})|^2 \, \, dt \bigg)^{2 p}, \]
            which is the first term in the required bound.
            
            We will now deal with the portion of the integral where $\Vert v_1 + \dots +  v_p - v_{p+1} - \dots - v_{2 p} \Vert < N^{-2/3}$. On this range, we will simply bound the minimum by $N$. Note that since we also have $\Vert v_i - u_i \Vert \le \frac{(\log N)^{100}}{N}$, we can see that \[\Vert u_1 + \dots + u_p - u_{p+1} - \dots - u_{2 p} \Vert < N^{-2/3} + \frac{2 p (\log N)^{100}}{N} \le \frac{2}{N^{2/3}}.\] Again bounding
            \[ \prod_{i=1}^{2p} \mathbbm{1}(\mathcal{A}_{u_i}) |F_K(e^{i u_i})| \mathbbm{1}(\mathcal{A}_{v_i}) |F_K(e^{i v_i})| \le \prod_{i=1}^{2p} \mathbbm{1}(\mathcal{A}_{u_i}) |F_K(e^{i u_i})|^2 + \prod_{i=1}^{2p} \mathbbm{1}(\mathcal{A}_{v_i}) |F_K(e^{i v_i})|^2 \]
            and using Proposition \ref{NTR2}, the contribution from this range is
            \begin{equation}\label{NearIntegerBound} \ll N \left( \frac{C \log \log N}{N} \right)^{2 p} \idotsint \prod_{i=1}^{2p} \mathbbm{1}(\mathcal{A}_{t_i}) |F_K(e^{i t_i})|^2 \, dt_1 \dots dt_{2 p}, \end{equation}
            where we are integrating over the set of $0 \le t_1, \dots, t_{2 p} \le 2 \pi$ satisfying $\Vert t_1 + \dots + t_p - t_{p+1} - \dots - t_{2 p} \Vert \le 2N^{-2/3}$. Unfortunately, it seems that some additional technical care is necessary in order to account for some wrap-around issues of variables around the circle. 
            
            We shall proceed to subdivide our range of integration. For a permutation $\sigma$ of length $2 p$, let 
            \[ I_{\sigma} := \big \{t_1, \dots, t_{2 p} : \Vert t_1 + \dots + t_p - t_{p+1} - \dots - t_{2 p}\Vert \le 2 N^{-2/3}, 0 < t_{\sigma(1)} < \dots < t_{\sigma(2 p)} < 2 \pi \big \}, \]
            so that these sets jointly cover our range of integration. We will sometimes identify $\sigma(2p+1)$ with $\sigma(1)$ and $\sigma(2p)$ with $\sigma(0)$. Moreover, for $1 \le h_1, \dots, h_{2 p} \le \log K$, let $I_{\sigma, h}$ be the subset of $I_{\sigma}$ where
            \[ e^{-h_i} \le \frac{\Vert t_{\sigma(i+1)} - t_{\sigma(i)} \Vert}{\pi} \le e^{-h_i+1} \]
            for all $i$, with the convention that the lower bound shall be omitted whenever $h_i = \log K$. Set 
            \[ \alpha(h) := \begin{cases}
            (\log N)^{12}, &\mbox{ if } 0 \le h < 0.01 \log K, \\
            (\log N)^{-2000}, &\mbox{ if } 0.01 \log K \le h \le \log K,
            \end{cases}\]
            and let $1 \le i_0 \le 2p$ be an index such that $h_{i_0}$ is minimal.
            
            Then by definition of $\mathcal{A}_t$ we can bound
            \begin{align*} \prod_{i = 1}^{2 p} \mathbbm{1}(\mathcal{A}_{t_i}) |F_K(e^{i t_i})|^2 &= \prod_{i = 1}^{2 p} \mathbbm{1}(\mathcal{A}_{t_{\sigma(i)}}) |F_K(e^{i t_{\sigma(i)}})|^2 \\ &\ll \bigg( \prod_{i=1}^{2 p} \mathbbm{1}(\mathcal{A}_{t_{\sigma(i)}}) \bigg) \exp \left( 2 \sum_{\substack{i = 1 \\ i \ne i_0}}^{2 p} \sum_{m=1}^{h_i} (Z_{t_{\sigma(i)}}(m) - 2 \sigma_m^2) + 4 \sum_{\substack{i = 1\\ i \ne i_0}}^{2 p} \sum_{m = 1}^{h_i} \sigma_m^2 + \right. \\  &\qquad \qquad \qquad \qquad \qquad \qquad  \left. + 2 \sum_{\substack{i = 1 \\ i \ne i_0}}^{2 p} \sum_{m = h_i + 1}^{\log K} Z_{t_{\sigma(i)}}(m) + 2 \sum_{m = 1}^{\log K} Z_{t_{\sigma(i_0)}}(m) \right) \\
            &\le \bigg( \prod_{\substack{i = 1 \\ i \ne i_0}}^{2 p} \alpha(h_i)^2 \exp \Big( 4 \sum_{m=1}^{h_i} \sigma_m^2 \Big) \! \bigg)  \exp \bigg( 2 \sum_{\substack{i = 1\\ i \ne i_0}}^{2 p} \sum_{m = h_i + 1}^{\log K}  \! \! \! Z_{t_{\sigma(i)}}(m) + 2 \sum_{m = 1}^{\log K}
            Z_{t_{\sigma(i_0)}}(m) \bigg) \\
            &\ll \bigg( \prod_{\substack{i = 1 \\ i \ne i_0 }}^{2 p} \alpha(h_i)^2 e^{2 h_i} \bigg) \exp \left( 2 \sum_{\substack{i = 1 \\ i \ne i_0}}^{2 p} \sum_{m = h_i + 1}^{\log K} Z_{t_{\sigma(i)}}(m) + 2 \sum_{m = 1}^{\log K} Z_{t_{\sigma(i_0)}}(m) \right).
            \end{align*}
            We note that, so far, we have not made any use of the parameter $h_{i_0}$. Denoting by $S_{n}$ the symmetric group on $n$ letters, we deduce that the expression in \eqref{NearIntegerBound} is
            \begin{align*} &\ll N \left( \frac{C \log \log N}{N} \right)^{2 p} \sum_{\sigma \in S_{2 p}} \sum_{0 \le h_1, \dots, h_{2 p} \le \log K} \bigg( \prod_{\substack{i = 1 \\ i \ne i_0}}^{2 p} \alpha(h_i)^2 e^{2 h_i} \bigg) \times \\ &\qquad \qquad \qquad \qquad \times \idotsint \limits_{I_{\sigma, h}} \exp \Bigg( 2 \sum_{\substack{i = 1 \\ i \ne i_0}}^{2 p} \sum_{m = h_i + 1}^{\log K} Z_{t_{\sigma(i)}}(m) + 2 \sum_{m = 1}^{\log K} Z_{t_{\sigma(i_0)}}(m) \Bigg) \, dt_1 \dots dt_{2 p}.  \end{align*}
            We will estimate the expectation of this quantity with the intent of later applying Markov's inequality. Note that by Proposition \ref{MomentsZ} the expectation of the expression inside the integral is
            \[ \ll \exp \Bigg( 2 \sum_{\substack{i = 1 \\ i \ne i_0 }}^{2 p} \sum_{m = h_i + 1}^{\log K} \sigma_m^2 + 2 \sum_{m = 1}^{\log K} \sigma_m^2 + 2 \sum_{1 \le i < j \le 2 p} \sum_{ e^{\max(h_i,h_j)} < k \le K } \frac{\cos(k(t_{\sigma(j)}- t_{\sigma(i)}))}{k} \Bigg) \]
            (with the convention that the term $k = 1$ is included in the last inner sum whenever $h_i = h_j = 0$). Note that if $i$ resp. $j$ is equal to $i_0$ then the last inner sum would start at $e^{h_j}$ resp. $e^{h_i}$, but since we chose $h_{i_0}$ to be minimal, we can include it in the maximum. Since the sum is empty whenever $h_i$ or $h_j$ are equal to $\log K$, Proposition \ref{NTR2} implies that this is in turn
            \begin{align}\label{TediousIndices} &\ll e^{O(p)} K \prod_{\substack{i = 1 \\ i \ne i_0}}^{2 p} \left( \frac{K}{e^{h_i}} \right) \exp \bigg( 6 \pi \sum_{\substack{ 1 \le i < j \le 2 p \\ h_i, h_j \ne \log K}} \frac{1}{\Vert t_{\sigma(j)} - t_{\sigma(i)} \Vert e^{\max(h_i,h_j)}} \bigg).
            \end{align}
            In order to understand the sum inside the exponential, we will split it into a sum over $i$ and $j$ such that $0 < t_{\sigma(j)} - t_{\sigma(i)} \le \pi$ and its complementary sum, where $\pi < t_{\sigma(j)} - t_{\sigma(i)} < 2 \pi$. Note that on the former range we have
            \begin{align*} \Vert t_{\sigma(j)} - t_{\sigma(i)} \Vert &= \Vert t_{\sigma(j)} - t_{\sigma(j-1)} \Vert + \dots + \Vert t_{\sigma(i+1)}-t_{\sigma(i)} \Vert \\ &\ge \pi (e^{-h_{j-1}} \mathbbm{1}(h_{j-1} \ne \log K) + \dots + e^{-h_{i+1}} \mathbbm{1}(h_{i+1} \ne \log K) + e^{-h_i}). \end{align*}
            On the latter range on the other hand, we can wrap around the other way to obtain that
            \begin{align*}
                \Vert t_{\sigma(j)} - t_{\sigma(i)} \Vert \ge \pi(e^{-h_{i-1}} \mathbbm{1}(h_{i-1} \ne \log K) + \dots + e^{-h_{j+1}} \mathbbm{1}(h_{j+1} \ne \log K) + e^{-h_j}),
            \end{align*}
            where the sum wraps from index $1$ around to index $2 p$. On the first range, we will first fix $i$ and sum over $j$, whereas on the second range we proceed conversely; since the arguments are symmetrical, we will restrict ourselves to the first case. We have
            \begin{align*} &\sum_{\substack{ 1 \le i < j \le 2 p \\ h_i, h_j \ne \log K \\ t_{\sigma(j)} - t_{\sigma(i)} \le \pi }} \frac{1}{\Vert t_{\sigma(j)} - t_{\sigma(i)} \Vert e^{\max(h_i,h_j)}} \\ &\qquad \quad \le \frac{1}{\pi} \sum_i \sum_{j} \frac{1}{\big(e^{-h_i} +e^{-h_{i+1}} \mathbbm{1}(h_{i+1} \ne \log K) + \dots + e^{-h_{j-1}} \mathbbm{1}(h_{j-1} \ne \log K) \big) e^{\max(h_i,h_j)}}, \end{align*}
            where the ranges of summation coincide. We will now subdivide the sum in the denominator into blocks of consecutive indices $l$ with $h_l \ne \log K$. We choose the lengths of these blocks inductively in the following way. We formally set $l_0 = i$. Supposing that the $m$-th block ends at an index $l_m$, we take $l_{m+1}$ (at most the largest value of $j$ in the above summation) minimal such that
            \[ e^{-h_{l_m+1}} \mathbbm{1}(h_{{l_m}+1} \ne \log K) + \dots + e^{-h_{l_{m+1}}} \mathbbm{1}(h_{l_{m+1}} \ne \log K) \ge e^{-h_i} \]
            if such index exists, and otherwise end the process. Note that within a given block, only the last index $l_m$ can have the property that $h_{l_m} \ge h_i$. We can therefore bound
            \[ \sum_{j \text{ in block } m} \frac{1}{e^{\max(h_i,h_j)}} \le 2 e^{-h_i}. \]
            Moreover, we have
            \[ e^{-h_i} +e^{-h_{i+1}} \mathbbm{1}(h_{i+1} \ne \log K) + \dots + e^{-h_{j-1}} \mathbbm{1}(h_{j-1} \ne \log K) \ge m e^{-h_i} \]
            uniformly over $j$ in block $m$. Since there are at most $2 p$ blocks, we therefore obtain, for any $i$, that
            \begin{align*} &\sum_{j} \frac{1}{\big(e^{-h_i} +e^{-h_{i+1}} \mathbbm{1}(h_{i+1} \ne \log K) + \dots + e^{-h_{j-1}} \mathbbm{1}(h_{j-1} \ne \log K) \big) e^{\max(h_i,h_j)}} \\ &\qquad \qquad \qquad \qquad \qquad \qquad \qquad \qquad \qquad \qquad \qquad \qquad \le \sum_{m = 1}^{2 p} \frac{2}{m} \le 2 \log(2p) + O(1). \end{align*}
            Summing this in turn over $i$ gives a bound $\le 4 p \log (2p) + O(p)$. Including the range where $\pi < t_{\sigma(j)} - t_{\sigma(i)} < 2 \pi$ and the factor $\frac{1}{\pi}$, this shows that the expression in \eqref{TediousIndices} is
            \[ \ll e^{O(p)} p^{48 p} K \prod_{\substack{i = 1 \\ i \ne i_0}}^{2 p} \left( \frac{K}{e^{h_i}} \right). \]
            So far, we have therefore shown that the expectation of the expression in \eqref{NearIntegerBound} is
            \begin{align*} &\ll N \left( C p^{24} \log \log N \right)^{2 p} \sum_{\sigma \in S_{2 p}} \sum_{0 \le h_1, \dots, h_{2 p} \le \log K} \bigg( \prod_{\substack{i = 1 \\ i \ne i_0}}^{2 p} \alpha(h_i)^2 e^{h_i} \bigg) \mathrm{vol}(I_{\sigma,h}).
            \end{align*}
            We can proceed in two ways in order to bound the volume of $I_{\sigma,h}$. One such way would be to simply note that after placing $t_{\sigma(i_0)}$ arbitrarily, our conditions on the distances of the $t$ variables gives
            \[ \mathrm{vol}(I_{\sigma,h}) \le e^{O(p)} \prod_{\substack{ i = 1 \\ i \ne i_0 }}^{2 p} e^{-h_i}. \]
            Another way would be to start in the same way, but for an arbitrary index $l$ with $l \ne i_0, i_0 + 1$ we instead use the condition $\Vert t_{1} + \dots + t_p - t_{p+1} - \dots - t_{2p} \Vert \le 2 N^{-2/3}$ on the variable $t_{\sigma(l)}$. This gives the bound
            \[ \mathrm{vol}(I_{\sigma,h}) \ll e^{O(p)} N^{-2/3} \min_{\substack{1 \le l \le 2p \\ l \ne i_0,i_0+1}} \prod_{\substack{i = 1 \\ i \ne i_0, \, l-1,\,  l }} \! \! \! \! e^{-h_i}. \]
            We hence obtain that the expectation of \eqref{NearIntegerBound} is
            \[ \ll N (C p^{24} \log \log N)^{2 p} \sum_{\sigma \in S_{2 p}} \sum_{0 \le h_1, \dots, h_{2 p} \le \log K} \bigg( \prod_{\substack{i = 1 \\ i \ne i_0}}^{2 p} \alpha(h_i)^2 \bigg) \min \Big \{ 1 , N^{-2/3} \min_{\substack{1 \le l \le 2p \\ l \ne i_0,i_0+1}} e^{h_l+h_{l+1}} \Big \}.  \]
            If there is $l \ne i_0, i_0+1$ such that $e^{h_l}, e^{h_{l+1}} \le N^{1/6}$, then we can bound the minimum in the above expression by $N^{-1/3}$ and bound the product over $\alpha$ trivially by $(\log N)^{48 p}$. If on the other hand there is no such $l$, then at least half of the remaining $2 p - 2$ indices must satisfy $e^{h_l} > N^{1/6}$ and thus in particular $h_l > 0.01 \log K$. In that case we bound the minimum trivially by $1$ and
            \[ \prod_{\substack{i = 1 \\ i \ne i_0}}^{2 p} \alpha(h_i)^2 \le (\log N)^{24p} (\log N)^{-4000 (p-1)} \le \left( (\log N)^{-988} \right)^{2 p}. \]
            We thus see, after noting that the number of permutations $\sigma$ is $\le e^{O(p)} p^{2p}$, that the expectation of \eqref{NearIntegerBound} is
            \[ \ll N^{2/3} (p^{25} (\log N)^{26})^{2 p} + N \left( \frac{p^{25}}{(\log N)^{986}} \right)^{2 p}. \]
            Therefore, Markov's inequality implies that with probability $\ge 1 - O((\log N)^{-2 p})$, the expression \eqref{NearIntegerBound} is
            \[ \ll N^{2/3} (p^{25} (\log N)^{27})^{2 p} + N \left( \frac{p^{25}}{(\log N)^{985}} \right)^{2 p}, \]
            thus implying the claim.
        \end{proof}
        
        \begin{proof}[Conclusion of the proof of Proposition \ref{CovarianceProposition}]
            Note first that, e.g. by \cite[Lemma 2.2]{SoundZaman1}, we have
            \[ \mathbb{E} \bigg[ \int_0^{2 \pi} |F_K(e^{i t})|^2 \, dt \bigg] \ll N, \]
            and so by Markov's inequality we have
            \[\int_0^{2 \pi} |F_K(e^{i t})|^2 \, dt \le N \log N \]
            with  probability $\ge 1 - O((\log N)^{-1}).$ As a consequence, we also have that the quantity on the left-hand side of Proposition \ref{TediousProposition} satisfies
            \begin{align*}
                &\max_{\frac{8 N}{7} < m \le \frac{4N}{3}} \sum_{\frac{8 N}{7} < n \le \frac{4N}{3}} \bigg| \iint \limits_{\substack{ [0, 2\pi] \times [0, 2 \pi] \\ \Vert t_2 - t_1 \Vert \le \frac{(\log N)^{100}}{N} }} \! \! \! \! \! \! \! \! \! \! \! \mathbbm{1}(\mathcal{A}_{t_1}) F_K(e^{i t_1}) \mathbbm{1}(\mathcal{A}_{t_2}) \overline{F_K(e^{i t_2})} e^{i t_1 (m+1) - i t_2 (n+1)} \! \! \! \! \sum_{N < k \le \frac{4N}{3}} \frac{e^{i(t_2-t_1) k}}{k} \, dt_1 \, dt_2 \bigg|^{2 p} \\
                &\qquad \qquad \qquad \qquad \qquad \qquad \qquad \qquad \qquad \qquad \ll N^{2/3} \left( p^{25} (\log N)^{27} \right)^{2 p} + N \left( \frac{p^{25}}{(\log N)^{985}} \right)^{2 p}
            \end{align*}
            with probability $\ge 1 - O((\log N)^{-1})$. Under this event, for any $\frac{8 N}{7} < m \le \frac{4 N}{3}$, we have 
            \[ \bigg| \iint \limits_{\substack{ [0, 2\pi] \times [0, 2 \pi] \\ \Vert t_2 - t_1 \Vert \le \frac{(\log N)^{100}}{N} }} \! \! \! \! \! \! \! \! \! \! \! \mathbbm{1}(\mathcal{A}_{t_1}) F_K(e^{i t_1}) \mathbbm{1}(\mathcal{A}_{t_2}) \overline{F_K(e^{i t_2})} e^{i t_1 (m+1) - i t_2 (n+1)} \! \! \! \! \sum_{N < k \le \frac{4N}{3}} \frac{e^{i(t_2-t_1) k}}{k} \, dt_1 \, dt_2 \bigg| \le \frac{1}{\log N}  \]
            for all but
            \[ \ll N^{2/3} \left( p^{25} (\log N)^{28} \right)^{2 p} + N \left( \frac{p^{25}}{(\log N)^{984}} \right)^{2 p}  \]
            values of $\frac{8 N}{7} < n \le \frac{4N}{3}$ and any integer $2 \le p \le N^{1/4}$. Taking $p = \left[ \frac{\log N}{B \log \log N} \right]$ for some absolute constant $B$, this is
            \[ \ll N^{2/3+106/B} + N^{1-1918/B}, \]
            which is certainly $\ll N^{7/10}$ by taking say $B = 6000$. Recalling \eqref{CovarianceRestrictA}, this gives the claim.
            
        \end{proof}

\bibliography{PM}

\newcommand{\etalchar}[1]{$^{#1}$}
\begin{thebibliography}{{Har}20a}

\bibitem[ASV{\etalchar{+}}21]{ASVZZ}
Daksh Aggarwal, Unique Subedi, William Verreault, Asif Zaman, and Chenghui
  Zheng.
\newblock A conjectural asymptotic formula for multiplicative chaos in number
  theory, 2021.

\bibitem[{Har}19]{Harper6}
Adam~J. {Harper}.
\newblock {On the partition function of the Riemann zeta function, and the
  Fyodorov--Hiary--Keating conjecture}.
\newblock {\em arXiv e-prints}, page arXiv:1906.05783, Jun 2019.

\bibitem[{Har}20a]{Harper8}
Adam~J. {Harper}.
\newblock {Almost sure large fluctuations of random multiplicative functions}.
\newblock {\em arXiv e-prints}, page arXiv:2012.15809, December 2020.

\bibitem[Har20b]{Harper1}
Adam~J. Harper.
\newblock Moments of random multiplicative functions, i: Low moments, better
  than squareroot cancellation, and critical multiplicative chaos.
\newblock {\em Forum of Mathematics, Pi}, 8:e1, 2020.

\bibitem[LTW13]{LTW1}
Yuk-Kam Lau, G\'{e}rald Tenenbaum, and Jie Wu.
\newblock On mean values of random multiplicative functions.
\newblock {\em Proc. Amer. Math. Soc.}, 141(2):409--420, 2013.

\bibitem[Mas21]{Mastrostefano1}
Daniele Mastrostefano.
\newblock An almost sure upper bound for random multiplicative functions on
  integers with a large prime factor, 2021.

\bibitem[SZ21]{SoundZaman1}
Kannan Soundararajan and Asif Zaman.
\newblock A model problem for multiplicative chaos in number theory, 2021.

\bibitem[You91]{RYoung}
Robert~M. Young.
\newblock {E}uler's {C}onstant.
\newblock {\em The Mathematical Gazette}, 75(472):187--190, 1991.

\end{thebibliography}
\bibliographystyle{alpha}

\end{document}